  \documentstyle[12pt]{article}
\newcommand{\numero}[1]{
\addtocounter{section}{1}
\begin{center}{\bf \thesection .\
#1\vspace{-.1in}}\end{center}
\setcounter{subsection}{0}
\setcounter{lemma}{0}\indent}

\newcommand{\subnumero}[1]{
\pagebreak[1]\begin{center}{\em #1}\nopagebreak\end{center}
}

\newcommand{\eop}{\hfill $/$\hspace*{-.1cm}$/$\hspace*{-.1cm}$/$\vspace{.1in}}

\newtheorem{lemma}{Lemma}[section]
\newtheorem{theorem}[lemma]{Theorem}

\newtheorem{corollary}[lemma]{Corollary}

\newtheorem{conjecture}{Conjecture}
\newtheorem{proposition}[lemma]{Proposition}

\newcommand{\mylabel}[1]{\label{#1}}

\newcommand{\zz}{{\bf Z}}

\newcommand{\Cc}{{\cal C}}
\newcommand{\Ee}{{\cal E}}

\newcommand{\Mm}{{\cal M}}

\newcommand{\Kk}{{\cal K}}

\newcommand{\tworightarrows}{\stackrel{\displaystyle \rightarrow}{\rightarrow}}

\begin{document}

\section*{Effective generalized Seifert-Van Kampen: how to calculate
$\Omega X$}

\noindent
Carlos Simpson\newline
CNRS, UMR 5580, Universit\'e Paul Sabatier, 31062 Toulouse CEDEX, France.

\noindent
email: carlos@picard.ups-tlse.fr

\bigskip

A central concept in algebraic topology since the 1970's has been that of {\em
delooping machine} \cite{BoardmanVogt} \cite{May} \cite{SegalTopology}.
Such a ``machine'' corresponds to a notion of $H$-space, or space with a
multiplication satisfying associativity, unity and inverse properties up to
homotopy in an appropriate way, including higher order coherences as first
investigated in \cite{Stasheff}.
A delooping machine is a specification of the extra homotopical structure
carried by the loop space $\Omega X$ of a connected basepointed topological
space $X$, exactly the structure allowing recovery of $X$ by a ``classifying
space'' construction.

The first level of structure is that the component set  $\pi _0(\Omega X)$ has a
structure of group  $\pi _1(X,x)$. Classically the {\em Seifert-Van
Kampen theorem} states that a pushout diagram of  connected spaces gives rise to
a pushout diagram of groups $\pi _1$. The  loop space construction $\Omega X$
with its delooping structure being the higher-order ``topologized''
generalization of $\pi _1$, an obvious question is whether a similar
Seifert-Van Kampen statement holds for  $\Omega X$.

The aim of this paper is to describe the operation underlying pushout of spaces
with loop space structure, answering the above question by giving a
Seifert-Van Kampen statement for delooping machinery. We work
with Segal's machine \cite{SegalInventiones}
\cite{Thomason}.  Our Seifert-Van Kampen statement is actually contained (in an
$n$-truncated version) in \cite{nCAT}. In the present paper, we don't
concentrate
on the formal aspects of this, but rather on the aspect of {\em effectivity}. It
turns out that the situation for higher homotopy is actually much better than
for $\pi _1$: one can effectively calculate the pushout of connected loop spaces
(of course in the nonconnected case, i.e. when the component groups are
nontrivial, one has the well-known effectivity problems for pushout of groups).
Again, rather than concentrate on formal aspects we present this in an applied
way as an algorithm to describe a finite cell complex representing the $n$-type
of $\Omega X$, for a given finite simply connected simplicial complex $X$.
Iterating in a relatively obvious way gives an algorithm for calculating $\pi
_i(X)$. At the end of the paper we briefly discuss the various formal aspects of
the situation and possible generalizations to other delooping machines.

The problem of giving an effective calculation of the $\pi
_i(X)$ of a simply
connected finite complex was first solved by E. Brown
(\cite{EBrown}, 1957)---and apparently also by A. Shapiro, unpublished.
Brown did
this by explicitly constructing the fibrations in the Postnikov tower of $X$. His
method is unrelated to generalized Seifert-Van Kampen.

In some sense, the notion of pushout and the Seifert-Van Kampen theorem
are well known and date back to J. H. C. Whitehead \cite{Whitehead} and D. Kan
\cite{Kan1} \cite{Kan2}. Whitehead showed in \cite{Whitehead}
that the pushout of classifying spaces for a diagram of groups, is a
classifying space for the pushout group, if both morphisms of groups are
injective. An immediate corollary is the same statement for simplicial groups.
This leads to the statement that the coproduct of simplicial groups corresponds
to coproduct of spaces, when at least one of the maps is a cofibration of
simplicial groups. Kan used this to describe, for any simplicial complex $X$, a
simplicial free group whose ``homology''  gives the $\pi _i(X)$ \cite{Kan1}
\cite{Kan2} \cite{Kan3}.  In fact, the realization of Kan's simplicial group is
equivalent to $\Omega X$ with its product structure \cite{Kan2}.

Kan's description seems to have a certain unconstructibility about it, since
one must do calculations in simplicial free groups, which doesn't seem to be a
finite process. In \cite{Curtis}, Curtis showed that the free groups can be
replaced by finite level lower central quotients (the $r$-th quotient of the
lower central series suffices to calculate $\pi _i$ for $i\leq \log _2 r$). This
makes calculations theoretically effective, and gave a second way  (after that
of \cite{EBrown}) to calculate the $\pi _i$. Note that it still doesn't give a
completely satisfactory answer for a model for $\Omega X$ since the component
groups are infinite. (For example it is  not entirely clear how one would go
about using Kan and Curtis to compute the cohomology ring of $\Omega
^2X$---whereas this becomes possible, in principle, with the algorithm we
present here as long as $X$ is $2$-connected.)

Recently Ellis \cite{Ellis} showed that if the $\pi _i$ are finite for
$i\leq n$ then the Curtis quotients can be further replaced by quotients which
are (nilpotent) finite groups, giving a simplicial finite group with the same
$n$-type as $X$. In this case of finite homotopy groups, one can say that the
Kan-Curtis method finally gives a theoretically reasonable model for $\Omega X$.

In the Whitehead-Kan-Curtis approach, the ``delooping structure'' which is
used is an actual associative group structure. In particular it is not
invariant under homotopy.  This cannot be seen as a default in and of itself.
However, experience in algebraic topology has shown that it is often useful to
work with the sole homotopy data which are involved in a loop space, rather than
an overly strictified but simpler structure. In this sense (although this is a
matter of opinion) there may be room for improvement in trying to do
Seifert-Van Kampen for homotopical delooping machinery rather than for
simplicial groups. For example,
Fiedorowicz noted using many examples in
\cite{Fiedorowicz}, that Whitehead's result doesn't generalize perfectly well
to the case of topological monoids.

Ronnie Brown has delved extensively into the question of pushouts and Van
Kampen theorems for $2$-types of spaces, or more generally for cases when the
loop space has a structure of ``crossed module'' (cf the list of references in
\cite{RBrown}). A related construction in the general case is Loday's
correspondence between cubes of spaces and $cat^n$-groups \cite{Loday}. Here,
some extra structure is required on the original space. In this case R. Brown
and Loday obtain a Seifert-Van Kampen theorem \cite{BrownLoday}.

A similar and essentially equivalent question is whether the classifying space
functor $B$ going in the other direction, from loop space data to space data, is
compatible with ``pushouts''. Whitehead \cite{Whitehead} proved that this is the
case, using the usual classifying space construction, for pushouts of discrete
groups along monomorphisms.  D. MacDuff in a special case  \cite{MacDuff} and
Fiedorowicz in general \cite{Fiedorowicz} prove a type of pushout theorem for a
pushout diagram of discrete monoids and the associated pushout diagram of
classifying spaces. There, there is an essential flatness assumption. The basic
reason for this is that the pushout in the sense of discrete monoids is not the
``right'' one.

One of the main approaches to calculation of loop spaces is based on James'
``reduced product construction'' \cite{James} which gives a way of calculating
$\Omega \Sigma X$ and---in generalizations using the delooping
machines---eventually leads to $\Omega ^n\Sigma ^nX$ for $1\leq n \leq \infty$.
This has led to many results on cohomology of loop spaces (see
\cite{CohenLadaMay} for example) but has the disadvantage that the suspension
tends to stabilize things.

Of course in rational homotopy theory, all the problems we discuss here have
been solved in a totally satisfactory way. It would go beyond our present scope
to get into references for that.

Beyond these
results, I couldn't find any treatment of the Seifert-Van Kampen question for
$\Omega X$ with its homotopical delooping structure, nor any explicit method of
calculating $\Omega X$.

\begin{center}
* \hspace*{1.5cm} *  \hspace*{1.5cm} *
\end{center}

This paper is based on the notion of {\em Segal category}. In its ``monoidal''
version (i.e. for categories with one object, which is what interests us in the
current paper) it is due to Segal  \cite{SegalTopology}
\cite{SegalInventiones}, and related notions are due to Stasheff
(\cite{Stasheff}, the first work on this matter), Boardman and Vogt
\cite{BoardmanVogt}, May \cite{May} and Thomason \cite{Thomason}
\cite{MayThomason} \cite{Thomason}. We will discuss (somewhat conjecturally) the
relationship between these notions at the end of the paper.

The generalization
to the case of several objects is immediate in Segal's point of view
(one must suppose that Segal was aware of this), and becomes possible for the
other delooping machines with the point of view adopted by Thomason and May in
\cite{MayThomason}, \cite{Thomason}.

In Boardman and Vogt
(\cite{BoardmanVogt}, p. 102) occurs a notion which is closely related to that
of Segal category, namely a simplicial set satisfying the ``restricted Kan
condition'', the Kan condition for horns obtained by deleting any but the first
or last faces. It seems likely that this notion is equivalent to the notion of
Segal category, but I haven't investigated this yet.
This occurence in Boardman and Vogt (1973) looks similar to the example in which
I personally came across the concept of Segal category in 1993, see
\cite{flexible}.

An extension of the linear-algebra version of Stasheff's notion
of $A_{\infty}$-algebra \cite{Stasheff}, to the case of several objects, is
given by Kontsevich in \cite{Kontsevich} where it is viewed as a weakened
version of the notion of ``differential graded category''.

In (Tamsamani \cite{Tamsamani}) the idea behind the notion of Segal category was
iterated to obtain the notion of $n$-category. The definitions of weak
$n$-category of Baez-Dolan \cite{BaezDolan} and Batanin \cite{Batanin} are
loosely based on May-type delooping machinery, although the connection doesn't
seem to be so direct as in Tamsamani's definition.

In \cite{nCAT} was proposed a ``generalized Seifert-Van Kampen theorem''
for the Poincar\'e $n$-groupoid $\Pi _n(X)$ defined by Tamsamani in
\cite{Tamsamani}.  The notion of (weak) $n$-groupoid is
quite analogous to Segal's notion of $n$-fold delooping machine---and
in fact the
point of view which we shall adopt below is that of Segal's delooping machine
rather than that of (weak) $n$-groupoid.
The key element in \cite{nCAT} was the operation $Cat$ allowing passage from an
$n$-precat to an $n$-category.  This operation is analogous to the description
of a group in terms of generators and relations. A natural question---stemming
from the apparently infinite nature of the operation $Cat$ as described in
\cite{nCAT}---is to what extent this can be made effective.

When viewed in an appropriate $n$-categorical sense (maybe for
$n=\infty$) any equivalence from spaces to objects with algebraic
structure, should satisfy Seifert-Van Kampen. Indeed, if the functor is an
equivalence then it automatically has to preserve coproducts. The real problem
is to find classes of algebraic objects and functors (equivalences) with the
property that one can actually compute the coproducts on the algebraic
side. Also, of course, the value of the functor on the algebraic side should
give more insight into the topology of the space than just the space itself.
It is in this light that we feel it interesting to dwell on effectivity.

In this  paper we will restrict to the question of effectively calculating the
$\pi _i(X)$ for simply connected finite simplicial complexes $X$.
This gives a point of view on the operation $Cat$ of \cite{nCAT} which has a
maximum of connectivity with homotopy theory. My main reason for looking at
this is to get some topological intuition for this operation. In particular, the
present paper could serve as an introduction to \cite{Tamsamani} and
\cite{nCAT}.

Our method gives yet another way (after E. Brown and Kan-Curtis) to
calculate the
$\pi _i$.  It doesn't seem to be any more ``realistically effective'' than the
previous methods.
The first thing which we obtain is an explicit (but large) finite complex
representing the $n$-type of $\Omega X$.
But even to calculate $H_i (\Omega X)$ the cell complex we
obtain has approximately
$$
n^{2n^2}
$$
times as many cells as $X$. This means that to calculate an
$H_{20}(\Omega X)$ would require more space-time than is available in the
entire universe!  Even worse, the well-known effectivity problems for $\pi _1$
crop in at every stage when we try to calculate $\pi _i$, and  in
the present state of the algorithm, the complexity is unbounded. It is
likely, of
course, that some improvements could be made. It is unclear whether this could
take things from a theoretical to an actually useful effectivity.

I would like to thank Andr\'e Hirschowitz, whose numerous questions and
suggestions led to the paper \cite{nCAT}; Zouhair Tamsamani whose work
\cite{Tamsamani} catalyzed the  thoughts
which go into the present paper; and Ronnie Brown who pointed out that a
generalized Van Kampen-type theorem is only interesting if it can be made
effective.

\numero{Segal categories}

The references for this section are \cite{Adams} \cite{SegalTopology} and
\cite{Thomason}.

Let $\Delta$ denote the simplicial category whose
objects are denoted $m$ for positive integers $m$, and where the morphisms
$p\rightarrow m$ are the (not-necessarily strictly) order-preserving maps
$$
\{ 0,1,\ldots , p\} \rightarrow \{ 0,1,\ldots , m\} .
$$
A morphism $1\rightarrow m$ sending $0$ to $i-1$ and $1$ to $i$ is called a
{\em principal edge} of $m$.  A morphism which is not injective is called a
{\em degeneracy}.

A {\em Segal precat} is a bisimplicial set
$$
A= \{ A_{p,k},
\;\; p,k\in \Delta \}
$$
(in other words a  functor $A: \Delta ^o\times \Delta ^o\rightarrow Sets$)
satisfying the {\em globular condition} that the simplicial set $k\mapsto
A_{0,k}$ is constant equal to a set which we denote by $A_0$ (called the set of
{\em objects}).

If $A$ is a Segal precat then for $p\geq 1$ we obtain a simplicial set
$$
k\mapsto A_{p,k}
$$
which we denote by $A_{p/}$.  This yields a simplicial collection of simplicial
sets.  One could instead look at {\em simplicial spaces} (i.e. take the
$A_{p/}$ to be spaces with $A_0$ discrete).  This gives an equivalent
theory, although there are degeneracy problems which apparently need to be
treated in an appendix in that case (\cite{MayThomason} \cite{Thomason}). We
note in passing that this necessary appendix is missing from \cite{Tamsamani}.

For each $m\geq 2$ there is a morphism of simplicial sets whose components are
given by the principal edges of $m$, which we call the {\em Segal map}:
$$
A_{m/}\rightarrow A_{1/} \times _{A_0} \ldots \times _{A_0} A_{1/}.
$$
The morphisms in the fiber product $A_{1/}\rightarrow A_0$ are alternatively the
inclusions $0\rightarrow 1$ sending $0$ to the object $1$, or to the object $0$.

We would like to think of the inverse image $A_{1/}(x,y)$ of a pair
$(x,y)\in A_0\times A_0$ by the two maps $A_{1/}\rightarrow A_0$ referred to
above, as the {\em simplicial set of maps from $x$ to $y$}.

We say that a Segal precat $A$ is  a {\em Segal category} if for all $m\geq
2$ the Segal maps
$$
A_{m/}\rightarrow A_{1/} \times _{A_0} \ldots \times _{A_0} A_{1/}.
$$
are weak equivalences of simplicial sets.

The main operation of this paper is a way of starting with a Segal precat
and enforcing the condition of becoming a Segal category, by forcing the
condition of weak equivalence on the Segal maps. As a general matter we will
call operations of this type $A\mapsto SeCat(A)$. We give an abstract general
discussion in \S 7 below, but for the body of the paper we do things
concretely.

Suppose $A$ is a Segal category. Then the simplicial set $p\mapsto \pi
_0(A_{p/})$ is the nerve of a category which we call $\tau _{\leq 1} A$.
We say that $A$ is a {\em Segal groupoid} if $\tau _{\leq 1} A$ is a groupoid.
This means that the $1$-morphisms of $A$ are invertible up to equivalence.

In fact we can make the same definition even for a Segal precat $A$: we define
$\tau _{\leq 1} A$ to be the simplicial set $p\mapsto \pi _0(A_{p/})$.

We can now describe exactly the situation envisaged in \cite{Adams}
\cite{SegalTopology}: a Segal category $A$ with only one object, $A_0 = \ast$.
We call this a {\em Segal monoid}. If $A$ is a groupoid then the homotopy
theorists' terminology is to say that it is {\em grouplike}.

\subnumero{Equivalences of Segal categories}

The basic intuition is to think of Segal categories as the natural weak version
of the notion of topological category. One of the main concepts in category
theory is that of a functor which is an ``equivalence of categories''.
This may be generalized to Segal categories. The same thing in the context
of $n$-categories is due to Tamsamani \cite{Tamsamani}.

We say
that a morphism $f:A\rightarrow B$ of Segal categories is an {\em equivalence}
if it is {\em fully faithful}, meaning that for $x,y\in A_0$ the map
$$
A_{1/}(x,y) \rightarrow B_{1/}(f(x),f(y))
$$
is a weak equivalence of simplicial sets; and  {\em essentially
surjective}, meaning that the induced functor of categories
$$
\tau _{\leq 1} (A)\rightarrow \tau _{\leq 1}(B)
$$
is surjective on isomorphism classes of objects. (Note that this induced
functor will be an equivalence of categories as a consequence of the
fully faithful condition.)

The homotopy theory that we are interested in is that of the category of Segal
categories modulo the above notion of equivalence. In particular, when we
search for the ``right answer'' to a question, it is only up to the above type
of equivalence. Of course when dealing with Segal categories having only one
object (as will actually be the case in what follows) then the essentially
surjective condition is vacuous and the fully faithful condition just amounts
to equivalence on the level of the``underlying space'' $A_{1/}$.

In order to have an appropriately reasonable point of view on the homotopy
theory of Segal categories one should look at the closed model structure
(discussed briefly in \S 7 below):  the right notion of weak morphism
from $A$ to $B$ is that of a morphism from $A$ to $B'$ where $B\hookrightarrow
B'$ is a fibrant replacement of $B$.  We don't want to get into this type
of question in the main part of the paper, since our aim is just to show how to
{\em calculate}  with these objects.

\subnumero{Segal's theorem}

We define the {\em realization} of a Segal category $A$ to be the space $|A|$
which is the realization of the bisimplicial set $A$.  Suppose $A_0 = \ast$.
Then we have a morphism
$$
|A_{1/} | \times {[}0,1] \rightarrow |A|
$$
giving a morphism
$$
|A_{1/}| \rightarrow \Omega |A|.
$$
The notation $|A_{1/}|$ means the realization of the simplicial set $A_{1/}$
and $\Omega |A|$ is the loop space based at the basepoint $\ast = A_0$.

\begin{theorem}
\mylabel{segal}
{\rm (G. Segal \cite{SegalTopology}, Proposition 1.5)}
Suppose $A$ is a Segal groupoid with one object. Then the morphism
$$
|A_{1/}| \rightarrow \Omega |A|.
$$
is a weak equivalence of spaces.
\end{theorem}

Refer to Segal's paper, or also May (\cite{MayFibs} 8.7), for a proof.

\subnumero{The translation with $n$-categories}

We mention briefly the relationship between the notions of Segal category and
$n$-category. This will not be used until \S 7 below. There, we will use it
to transfer some results on $n$-categories to results on Segal categories
(transfering the proof techniques, allowing us to skip the proofs).  On the
other
hand, the reader should also be able, via this translation, to use the present
paper as an introduction to \cite{Tamsamani}, \cite{nCAT}.

Tamsamani's definition of $n$-category is recursive. The basic idea is to use
the same definition as above for Segal category, but where the $A_{p/}$ are
themselves $n-1$-categories. The appropriate condition on the Segal maps is
the condition of equivalence of $n-1$-categories, which in turn is defined
(inductively) in the same way as the notion of equivalence of Segal categories
explained above.

Tamsamani shows
\footnote{
Actually the proof in \cite{Tamsamani} using simplicial spaces, is missing a
discussion of ``whiskering'' as is standard in delooping and classifying space
constructions (cf \cite{SegalTopology} \cite{May}, \cite{MayThomason},
\cite{Thomason}). Alternatively the proof works as it is if ``spaces'' are
replaced by ``simplicial sets''.} that the homotopy category of $n$-groupoids is
the same as that of $n$-truncated spaces. The two relevant functors are the
realization and Poincar\'e $n$-groupoid $\Pi _n$ functors. Applying this to
the $n-1$-categories $A_{p/}$ we obtain the following relationship. An
$n$-category $A$ is said to be {\em $1$-groupic} (notation introduced in
\cite{limits}) if the $A_{p/}$ are  $n-1$-groupoids. In this case, replacing
the $A_{p/}$ by their realizations $|A_{p/}|$ we obtain a simplicial space
which satisfies the Segal condition. Conversely if $A_{p/}$ are spaces or
simplicial sets then replacing them by their $\Pi _{n-1}(A_{p/})$ we obtain
a simplicial collection of $n-1$-categories, again satisfying the Segal
condition.  These constructions are not quite inverses because
$$
| \Pi _{n-1}(A_{p/})| = \tau _{\leq n-1} (A_{p/})
$$
is the Postnikov truncation. If we think (heuristically) of setting
$n=\infty$ then we get inverse constructions. Thus---in a sense which I will
not currently make more precise than the above discussion---one can say that
Segal categories are the same thing as $1$-groupic $\infty$-categories.

The passage from simplicial sets to Segal categories is the same as the
inductive passage from $n-1$-categories to $n$-categories. In \cite{nCAT}
was introduced the notion of {\em $n$-precat}, the analogue of the above
Segal precat.  Noticing that the results and arguments in \cite{nCAT} are
basically organized into one gigantic inductive step passing from $n-1$-precats
to $n$-precats, the same step applied only once works to give the analogous
results in the passage from simplicial sets to Segal precats.

The notion of Segal category thus presents, from a technical point of view,
an aspect of a ``baby'' version of the notion of $n$-category.  On the other
hand, it allows a first introduction of homotopy going all the way up to
$\infty$ (i.e. it allows us to avoid the $n$-truncation inherent in the
notion of $n$-category).

One can easily imagine combining the two into a notion of ``Segal
$n$-category'' which would be an $n$-simplicial simplicial set satisfying the
globular condition at each stage. It is interesting and historically
important to note that the notion of Segal $n$-category with only one
$i$-morphism for each $i\leq n$, is the same thing as the notion of {\em
$n$-fold delooping machine}.
This translation comes out of Dunn \cite{Dunn}, which apparently dates
essentially back to 1984.  In retrospect it is not too hard to see
how to go from Dunn's notion of $E_n$-machine, to Tamsamani's notion of
$n$-category, simply by relaxing the conditions of having only one object.
Metaphorically, $n$-fold delooping machines correspond to the
Whitehead tower, whereas $n$-groupoids correspond to the Postnikov tower.

\numero{How to calculate $\Omega X$}

Suppose $X$ is a simplicial set with $X_0 = X_1 = \ast$, and with finitely
many nondegenerate simplices.  Fix $n$. We will obtain, by iterating an
operation
closely related to the operation $Cat$ of \cite{nCAT}, a finite complex
representing the $n$-type of $\Omega X$. See \S 7 for a more abstract
version of this operation which we denote $SeCat$.

Suppose $A$ is a Segal precat with $A_0 = \ast$. We say that $A$ is {\em
$(m,k)$-arranged} if the Segal map
$$
A_{m/} \rightarrow A_{1/} \times  \ldots \times A_{1/}
$$
induces isomorphisms on $\pi _i$ for $i< k$ and a surjection on $\pi _k$.

Note that for $l\geq k$, adding $l$-cells to $A_{m/}$ or $l+1$-cells to
$A_{1/}$ doesn't affect this property.

\begin{theorem}
\mylabel{bound1}
If $A$ is a Segal precat with $A_0= \ast$ and $A_{1/}$ connected, such that
$A$ is $(m,k)$-arranged for all $m+k \leq n$ then there exists a morphism
$A\rightarrow A'$ such that:
\newline
(1)\,\, the morphism $|A|\rightarrow |A'|$ is a weak equivalence;
\newline
(2)\,\, $A'$ is a Segal
groupoid; and
\newline
(3)\,\, the map of simplicial sets $A_{m/} \rightarrow A'_{m/}$ induces an
isomorphism on $\pi _i$ for $i+m < n$.
\end{theorem}

The proof of this theorem  will be given in \S 6 below.

\begin{corollary}
\mylabel{bound2}
Suppose $A$ is a Segal precat with $A_0= \ast$ and $A_{1/}$ connected, such that
$A$ is $(m,k)$-arranged for all $m+k \leq n$. Then the natural morphism
$$
|A_{1/}| \rightarrow \Omega | A|
$$
induces an isomorphism on $\pi _i$ for $i < n-1$.
\end{corollary}
{\em Proof:}
Use Theorem \ref{bound1} to obtain a morphism $A\rightarrow A'$ with the
properties stated there (which we refer to as (1)--(3)).  We have a diagram
$$
\begin{array}{ccc}
|A_{1/}| & \rightarrow & \Omega | A| \\
\downarrow & & \downarrow \\
|A'_{1/}| & \rightarrow & \Omega | A'|
\end{array} .
$$
By property (1) the vertical morphism on the right is a weak equivalence.
By property (2) and Theorem \ref{segal} the morphism on the bottom is a weak
equivalence. By property (3) the vertical morphism on the right induces
isomorphisms on $\pi _i$ for $i<n-1$. This gives the required statement.
\eop

In view of Corollary \ref{bound2}, in order to calculate the $n$-type of
$\Omega |A|$ we just have to change $A$ by pushouts preserving the weak
equivalence type of $|A|$ in such a way that $A$ is $(m,k)$-arranged for all
$m+k \leq n+2$.

{\em Remark:}
Theorem \ref{bound1} is used only in order to show that our procedure
actually gives the right answer. In particular it doesn't need an effective
proof, and we will make free use of infinite sequences of operations, during the
proof in \S 6 below. It is for this reason that it seemed like a good idea to
put that proof in  a separate section below, since it would seem out of place in
our current ``effective'' world.

We now define an operation where we try to ``arrange'' $A$ in degree $m$.
This operation is inspired by the operation $Raj$ of \cite{nCAT}. We
call this
$$
A \mapsto Arr (A, m).
$$
Fix $m$ in what follows.
Let $C$ be the mapping cone of the Segal map
$$
A_{m/} \rightarrow A_{1/}\times \ldots \times A_{1/}.
$$
To be precise, as a bisimplicial set
$$
C = (I \times A_{m/}) \cup ^{\{ 1\} \times A_{m/}} (A_{1/}\times \ldots \times
A_{1/}),
$$
where $I$ is the standard simplicial interval, and the notation is coproduct
of bisimplicial sets (note also that the globular condition is preserved, so it
is a coproduct of Segal precats). Note that $\{ 1 \} \times A_{m/}$ denotes the
second endpoint of the interval crossed with $A_{m/}$. We have morphisms
$$
A_{m/} \stackrel{a}{\hookrightarrow} C
\stackrel{b}{\rightarrow} A_{1/}\times \ldots \times
A_{1/}, $$
the morphism $a$ being the inclusion of $\{ 0\} \times A_{m/}$ into
$I \times A_{m/}$ (thus it is a cofibration i.e. injection of simplicial
sets) and the second morphism $b$ coming from the projection  $I\times A_{m/}
\rightarrow A_{m/}$. The second morphism $b$ is a weak equivalence.

We now define $Arr (A, m)$ as follows. For any $p$, let
$$
Arr (A,m)_{p/} := A_{p/} \cup ^{(\bigcup A_{m/})} \left( \bigcup _{p\rightarrow
m}
 C \right)
$$
be the combined coproduct of $A_{p/}$ with several copies of the
morphism $a:A_{m/} \rightarrow C$ , one copy for each map $p\rightarrow m$ not
factoring through a principal edge (see below for further discussion of this
condition), these maps inducing $A_{m/} \rightarrow A_{p/}$.

We need to define $Arr (A,m)$ as a Segal precat, i.e. as a bisimplicial set.
For this we need morphisms of functoriality
$$
Arr (A, m)_{p/} \rightarrow Arr (A, m)_{q/}
$$
for any $q\rightarrow p$. These are defined as follows.
We consider a component of $Arr (A, m)_{p/}$ which is a copy of $C$
attached along a map $A_{m/}\rightarrow A_{p/}$ corresponding to $p\rightarrow
m$ which doesn't factor through a principal edge.  If the composed map
$q\rightarrow p\rightarrow m$ doesn't factor through a principal edge
then the component $C$ maps to the corresponding component of
$Arr (A, m)_{1/}$. If the map does factor through a principal edge
$q\rightarrow 1\rightarrow m$ then we obtain a map
$C \rightarrow A_{1/}$ (the component of the map $b$ corresponding to this
principal edge). Compose with the map $A_{1/}\rightarrow A_{q/}$ to obtain a
map $C\rightarrow A_{q/}$. Note that if the map further factors
$$
q\rightarrow 0 \rightarrow 1 \rightarrow m
$$
then the map $A_{1/}\rightarrow A_{q/}$ factors through the basepoint
$$
A_{1/}\rightarrow A_0 \rightarrow A_{q/},
$$
and our map on $C$ factors through the basepoint. This factorization doesn't
depend on choice of principal edge containing the map $0\rightarrow m$.

One can verify that this prescription defines a functor $p\mapsto Arr
(A,m)_{p/}$ from $\Delta$ to simplicial sets. This verification will be a
consequence of the more conceptual description which follows.

Let $h(m)$ denote the simplicial set representing the standard $m$-simplex;
it is the contravariant functor on $\Delta$ represented by the object $m$.
Let $\Sigma (m)\subset h(m)$ be the subcomplex which is the union of the
principal edges.

\noindent
{\bf Notation:} If $X$ is a simplicial set and $B$ is another simplicial set
denote by $X\otimes B$ the bisimplicial set exterior product, defined by
$$
(X\otimes B)_{p,q} := X_p \times B_q.
$$

If $B$ is any simplicial set then putting $h(m)$ or $\Sigma (m)$ in the first
variable, we obtain an inclusion of bisimplicial sets which we denote
$$
\Sigma (m)\otimes B \hookrightarrow h(m)\otimes B .
$$
Note that these bisimplicial sets are not Segal precats because they don't
satisfy the globular condition (they are not constant over $0$ in the first
variable). However, that the morphism of simplicial sets
$$
(\Sigma (m)\otimes B)_{0/}\hookrightarrow  (h(m)\otimes B)_{0/}
$$
is an isomorphism because $\Sigma$ contains all of the vertices.

If $A$ is a Segal precat then a morphism $h(m)\otimes B
\rightarrow A$ is the same thing as a morphism $B\rightarrow A_{m/}$. Similarly,
a morphism $$
\Sigma (m)\otimes B \rightarrow A
$$
is the same thing as a morphism
$$
B\rightarrow A_{1/} \times _{A_0} \ldots \times _{A_0}A_{1/}.
$$

The morphism of realizations
$$
|\Sigma (m)\otimes B | \rightarrow | h(m) \otimes B|
$$
is a weak equivalence. To see this note that it is the product of $|B|$ and
$$
|\Sigma (m)| \rightarrow | h(m) |,
$$
and this last morphism is a weak equivalence (it is the inclusion from the
``spine'' of the $m$-simplex to the $m$-simplex; both are contractible).

Suppose $B'\subset B$ is an injection of simplicial sets. Put
$$
U:= \left( \Sigma (m)\otimes B\right) \cup ^{\Sigma (m) \otimes B'}
\left( h(m)\otimes
B'\right) ,
$$
and
$$
V:= h(m)\otimes B.
$$
We have an injection $U\hookrightarrow V$.
If $A$ is a Segal precat then a map $U\rightarrow A$ consists of a
commutative diagram
$$
\begin{array}{ccc}
B' & \rightarrow & B \\
\downarrow && \downarrow \\
A_{m/} & \rightarrow & A_{1/} \times _{A_0}\ldots \times _{A_0} A_{1/} .
\end{array}
$$
The inclusion
$$
\Sigma (m) \otimes B \rightarrow U
$$
induces a weak equivalence of realizations, because of the fact that the
inclusion
$\Sigma (m) \otimes B'\rightarrow h(m)\otimes B'$ does. Therefore the morphism
$|U|\rightarrow |V|$ is a weak equivalence.

We can now interpret our operation $Arr (A, m)$ in these terms. Applying the
previous paragraph to the inclusion $A_{m/} \hookrightarrow C$, we obtain
an inclusion of bisimplicial sets
$U\hookrightarrow V$. We get a map $U\rightarrow A$ corresponding to the diagram
$$
\begin{array}{ccc}
A_{m/} & \rightarrow & C \\
\downarrow && \downarrow \\
A_{m/} & \rightarrow & A_{1/} \times _{A_0}\ldots \times _{A_0} A_{1/} .
\end{array}
$$
The left vertical arrow is the identity map, the top arrow is $a$ and the
right vertical arrow is $b$. The bottom arrow is the Segal map.

It is easy to see that
$$
Arr (A,m) = A \cup ^U V.
$$
In passing, this proves associativity of the previous formulas for
functoriality of $Arr (A,m)$.

We get
$$
|Arr (A,m) |= |A| \cup ^{|U|} |V|.
$$
Since $|U|\rightarrow |V|$ is a weak equivalence, this implies the

\begin{lemma}
\mylabel{triviality}
The morphism induced by the above inclusion on realizations,
$$
|A| \hookrightarrow |Arr (A,m)|
$$
is a weak equivalence of spaces.
\end{lemma}
\eop

The key observation is the following proposition.

\begin{proposition}
\mylabel{arrangement}
Suppose $A$ (with $A_0 = \ast$ and $A_{1/}$ connected) is
$(m, k-1)$-arranged and $(p,k)$-arranged for some
$p\neq m$. Then  $Arr (A, m)$ is $(p,k)$-arranged and $(m,k)$-arranged.
\end{proposition}
{\em Proof:}
Keep the hypotheses of the proposition.  Let $C$ be the cone occuring in the
construction $Arr (A, m)$. Denote $B:= Arr (A, m)$. Then the map
$$
a:A_{m/} \hookrightarrow C
$$
is weakly equivalent to a map obtained by adding cells of dimension $\geq k$ to
$A_{m/}$.  This is by the condition that $A$ is $(m, k-1)$-arranged. Let
$h_1,\ldots , h_u$ be the $k$-cells that are attached to $A_{m/}$ to give $C$.

We first show that $B$ is $(p,k)$-arranged.
Note that $B_{p/}$ is obtained from $A_{p/}$ by attaching a certain number of
$k$-cells, $h^{p\rightarrow m}_i$ for $i=1,\ldots , u$ indexed by the maps
$p\rightarrow m$ not factoring through the principal edges of $m$; plus some
cells of dimension $\geq k+1$. The higher-dimensional  cells don't have any
effect on the question of whether $A$ is $(p,k)$-arranged.

On the other hand, $B_{1/}$ is obtained from $A_{1/}$ by attaching cells
$h^{1\rightarrow m}_i$ for $i=1,\ldots , u$ and indexed by the maps
$1\rightarrow m$ not factoring through the principal edges. Note here that
these maps cannot be degenerate, thus they are the non-principal edges.

Now $B_{1/} \times \ldots \times B_{1/}$ (product of $p$-copies) is obtained
from  $A_{1/} \times \ldots \times A_{1/}$ by adding $k$-cells indexed as
$\nu _{1\rightarrow p}(h^{1\rightarrow m}_i)$ where the indexing
$1\rightarrow p$ are principal edges and $1\rightarrow m$ are non-principal
edges. Then by adding cells of dimension $\geq k+1$ which have no effect on the
question. The notation $\nu _{1\rightarrow p}$ refers to the map
$$
A_{1/} \rightarrow A_{1/}\times \ldots \times A_{1/}
$$
putting the base point (i.e. the degeneracy of the unique point in $A_0$) in all
of the factors except the one corresponding to the map $1\rightarrow p$.

For every principal edge $1\rightarrow p$ there is a unique degeneracy
$p\rightarrow 1$ inducing an isomorphism $1\rightarrow 1$ and this establishes
a bijection between principal edges and degeneracies. Thus we may rewrite our
indexing of the $k$-cells attached to the product above as $\nu _{p\rightarrow
1}(h^{1\rightarrow m}_i)$.

Now for every pair $(p\rightarrow 1, 1\rightarrow m)$
the composition $p\rightarrow m$ is a degenerate morphism, not factoring
through a principal edge; and these degenerate morphisms are all different for
different pairs $(p\rightarrow 1, 1\rightarrow m)$. Thus $B_{p/}$ contains a
$k$-cell $h^{p\rightarrow m}_i$ for each $i=1,\ldots , u$ and each of these
maps $p\rightarrow m$ (plus possibly other cells for other maps $p\rightarrow
m$ but we don't use these).  Take such a cell $h^{p\rightarrow m}_i$, and look
at its image in $B_{1/} \times \ldots \times B_{1/}$ by the Segal map.
The projection to any factor $1\rightarrow p$ other than the one which splits
the degeneracy $p\rightarrow 1$, is totally degenerate coming from a
factorization $1\rightarrow 0\rightarrow 1$, hence goes to the unique basepoint.
The projection to the unique factor which splits the degeneracy is just the
cell $h^{1\rightarrow m}_i$. Thus the projection of our cell
$h^{p\rightarrow m}_i$ to the product is exactly the cell
$\nu _{p\rightarrow
1}(h^{1\rightarrow m}_i)$.   This shows that all of the new $k$-cells which
have been added to the product, are lifted as new $k$-cells in $B_{p/}$.
Together with the fact that $A_{p/} \rightarrow A_{1/} \times \ldots \times
A_{1/}$ was an isomorphism on $\pi _i$ for $i< k$ and a surjection for $i=k$,
we obtain the same property for $B_{p/} \rightarrow B_{1/}\times \ldots \times
B_{1/}$.  Note that the further $k$-cells which are attached to $B_{p/}$ by
morphisms $p\rightarrow m$ other than those we have considered above, don't
affect this property. (In general, attaching $k$-cells to the domain of a map
doesn't affect this property, but attaching cells to the range can affect it,
which was why we had to look carefully at the cells attached to $B_{1/}$).
This completes the proof that $B$ remains $(p,k)$-arranged.

We now prove that $B$ becomes $(m,k)$-arranged.
Note that $B_{m/}$ is obtained by first adding on $C$ to $A_{m/}$ via the
identity map $m\rightarrow m$; then adding some other stuff which we treat in a
minute. The Segal map for $B$ maps this copy of $C$ directly into
$A_{1/} \times \ldots A_{1/}$. The fact that $C$ is a
mapping cone for the Segal map means that the map
$$
C\rightarrow A_{1/} \times \ldots A_{1/}
$$
is a homotopy equivalence. In particular, it is bijective on $\pi _i$ for $i<k$
and surjective for $i=k$.

Now $B_{m/}$ is obtained from $C$ by adding various cells to $C$ along
degenerate
maps $m\rightarrow m$. The new $k$-cells which are added to $B_{1/}
\times \ldots \times B_{1/}$ ($m$-factors this time) are lifted to cells in
$B_{m/}$ added to $C$ via the degeneracies $m\rightarrow m$ which factor
through a principal edge. The argument is the same as above and we don't repeat
it. We obtain that $B$ is $(m,k)$-arranged.

This completes the proof of the proposition.
\eop

\begin{corollary}
\mylabel{thinkIcan}
Fix $n$, and suppose $A$ is a Segal precat with $A_0 = \ast$ and $A_{1/}$
connected. By applying the operations $A \mapsto Arr (A, m)$ for various $m$, a
finite number of times (less than $(n+2)^2$) in a predetermined way,
we can effectively get to a morphism of Segal precats $A\rightarrow B$ such
that
$$
|A| \rightarrow |B|
$$
is a weak equivalence of spaces, and such that $B$ is $(m,k)$-arranged for all
$m+k \leq n+2$.
\end{corollary}
{\em Proof:}
By Corollary \ref{triviality} any successive application of the operations
$A\mapsto Arr (A, m)$ yields a morphism $|A|\rightarrow |B|$ which is a weak
equivalence of spaces. By Proposition \ref{arrangement} it suffices, for
example,
to successively apply $Arr (A,i)$ for $i= 2, 3,\ldots , n+2$ and to repeat this
$n+2$ times.
\eop

\begin{corollary}
\mylabel{OmegaDone}
Fix $n$, and suppose $A$ is a Segal precat with $A_0 = \ast$ and $A_{1/}$
connected. Let $B$ be the result of the operations of Corollary
\ref{thinkIcan}. Then the $n$-type of the simplicial set $B_{1/}$ is equivalent
to the $n$-type of $\Omega |A |$.
\end{corollary}
{\em Proof:}
Apply Corollaries \ref{bound2} and \ref{thinkIcan}.
\eop

Going back to the original situation, suppose $X$ is a simplicial set
with finitely many nondegenerate
simplices, with $X_0=X_1 = \ast$. Apply the above letting $A$ be $X$
considered as a Segal precat constant in the second variable, in other words
$$
A_{p,k} := X_p.
$$

\begin{corollary}
\mylabel{OmegaDone2}
Fix $n$, and suppose $X$ is a simplicial set with finitely many nondegenerate
simplices, with $X_0=X_1 = \ast$. Let $A$ be $X$
considered as a Segal precat. Let $B$ be the result of the operations of
Corollary \ref{thinkIcan}. Then the $n$-type of the simplicial set $B_{1/}$ is
equivalent to the $n$-type of $\Omega X$.
 \end{corollary}
{\em Proof:}
An immediate restatement of \ref{OmegaDone}.
\eop

{\em Remark:} Any finite region of the Segal precat $B$ is effectively
computable. In fact it is just an iteration of operations pushout and mapping
cone, arranged in a way which depends on combinatorics of simplicial sets.
Thus the $n+1$-skeleton of the simplicial set $B_{1/}$ is effectively
calculable (in fact, one could bound the number of simplices in $B_{1/}$).

\begin{corollary}
Fix $n$, and suppose $X$ is a simplicial set with finitely many
nondegenerate simplices, with $X_0=X_1 = \ast$. Then we
can effectively calculate $H_i (\Omega |X|, \zz )$ for $i\leq n$.
\end{corollary}
{\em Proof:}
Immediate from above.
\eop

In some sense this corollary is the ``most effective'' part of the present
paper, since we can get at the calculation after a bounded number of
easy steps of the form $A\mapsto Arr (A,m)$.

\numero{Getting $A_{1/}$ to be connected}

We treat here the question of how to arrange $A$ on the level of $\tau
_{\leq 1}(A)$.

We define operations $Arr ^{0\, {\rm only}}(A,m)$ and
$Arr ^{1\, {\rm only}}(A,m)$. These consist of doing the operation $Arr (A,m)$
but instead of using the entire  mapping cone $C$, only adding on $0$-cells
to $A_{m/}$ to get a surjection on $\pi _0$; or only adding on $1$-cells to get
an injection on $\pi _0$.  Note in the second case that we {\em don't} add
extra $0$-cells. This is an important point, because if we added further
$0$-cells every time we added some $1$-cells, the process would never stop.

To define  $Arr ^{0\, {\rm only}}(A,m)$, use the same construction as for $Arr
(A,M)$ but instead of setting $C$ to be the mapping cone, we put
$$
C':= A_{m/} \cup sk_0(A_{1/} \times \ldots \times A_{1/}.
$$
Here $sk_0$ denotes the $0$-skeleton of the simplicial set, and $im$ means the
image under the Segal map. Let $C\subset C'$ be a subset where we choose
only one point for each  connected component of the product.
With this $C$ the
same construction as previously gives $Arr ^{0\, {\rm only}}(A,m)$.

With the subset $C\subset C'$ chosen as above (note that this choice can
effectively be made) the resulting simplicial set
$$
p\mapsto \pi _0\left( Arr ^{0\, {\rm only}}(A,m)_{p/} \right)
$$
may be described only in terms of the simplicial set
$$
p\mapsto \pi _0(A_{p/}).
$$
That is to say, this operation $Arr ^{0\, {\rm only}}(A,m)$ commutes with the
operation of componentwise applying $\pi _0$. We formalize this as
$$
\tau _{\leq 1}Arr ^{0\, {\rm only}}(A,m)
=\tau _{\leq 1}Arr ^{0\, {\rm only}}(\tau _{\leq 1}A,m).
$$

To define $Arr ^{1\, {\rm only}}(A,m)$, let $C$ be the cone of the map from
$A_{m/}$ to
$$
im(A_{m/}) \cup sk_1(A_{1/} \times \ldots \times A_{1/})^o
$$
where  $sk_1(A_{1/} \times \ldots \times A_{1/})^o$ denotes the union of
connected components of the $1$-skeleton of the product, which touch
$im(A_{m/})$. In this case, note that the inclusion
$$
A_{m/} \hookrightarrow C
$$
is $0$-connected (all connected components of $C$ contain elements of $A_{m/}$).
Using this $C$ we obtain the operation $Arr ^{1\, {\rm only}}(A,m)$.
It doesn't introduce any new connected components in the
new simplicial sets $A'_{p/}$, but may connect together some components which
were disjoint in $A_{p/}$.

Again, the operation $Arr ^{1\, {\rm only}}(A,m)$ commutes with truncation:
we have
$$
\tau _{\leq 1}Arr ^{1\, {\rm only}}(A,m)
=\tau _{\leq 1}Arr ^{1\, {\rm only}}(\tau _{\leq 1}A,m).
$$

Our goal in this section is to find a sequence of operations which makes $\tau
_{\leq 1}(A)_1$ become trivial (equal to $\ast$).  In view of this, and the
above commutations, we may henceforth work with simplicial sets (which we
denote $U= \tau _{\leq 1}A$ for example) and use the above operations
followed by
the truncation $\tau _{\leq 1}$ as modifications of the simplicial set $U$.
We try to obtain $U_1=\ast$.  This corresponds to making $A_{1/}$  connected.

Our operations have the following interpretation. The operation
$$
U\mapsto \tau _{\leq 1}Arr ^{0\, {\rm only}}(U,2)
$$
has the effect of formally adding to $U_1$ all binary products of pairs of
elements in $U_1$.  (We say that a binary product of $u,v\in U_1$ is defined if
there is an element $c\in U_2$ with principal edges $u$ and $v$ in $U_1$; the
product is then the image $w$ of the third edge of $c$).

The operation
$$
U\mapsto \tau _{\leq 1}Arr ^{1\, {\rm only}}(U,2)
$$
has the effect of identifying $w$ and $w'$ any time both $w$ and $w'$ are
binary products of the same elements $u,v$.

The operation
$$
U\mapsto \tau _{\leq 1}Arr ^{0\, {\rm only}}(U,3)
$$
has the effect of introducing, for each triple $(u,v,w)$, the
various binary products one can make (keeping the same order) and giving  a
formula
$$
(uv)w = u(vw)
$$
for certain of the binary products thus introduced.

It is somewhat unclear whether blindly applying the composed operation
$$
U\mapsto \tau _{\leq 1}
Arr ^{1\, {\rm only}}(\tau _{\leq 1}Arr ^{0\, {\rm only}}(U,3),2)
$$
many times must automatically lead to $U_1=\ast$ in case the actual fundamental
group is trivial.   This is because in the process of adding
the associativity, we also add in some new binary products; to which
associativity might then have to be applied in order to get something trivial,
and so on.

If the above doesn't work, then we may need a slightly revised version of the
operation $Arr ^{0\, {\rm only}}(U,3)$ where we add in only certain triples
$u,v,w$.  This can be accomplished by choosing a subset of the original $C$ at
each time.  Similarly for the
$Arr ^{0\, {\rm only}}(U,2)$ for binary products.
We now obtain a situation where
we have operations which effect the appropriate changes on $U$ corresponding to
all of the various possible steps in an elementary proof that the associative
unitary monoid generated by generators $U_1$ with relations $U_2$, is trivial.
Thus if we have an elementary proof that the associative unitary monoid
generated by $U_1$ with relations $U_2$ is trivial, then we can read off
from the
steps in the proof, the necessary sequence of operations to apply to get
$U_1=\ast$.  On the level of $A$ these same steps will result in a new $A$ with
$A_{1/}$ connected.

In our case we are interested in the {\em group completion} of the monoid: we
want to obtain the condition of being a Segal groupoid not just a Segal
category. It is possible that the simplicial set $X$ we start with would yield
a monoid which is not a group, when the above operations are applied. To fix
this, we take note of another operation which can be applied to $A$ which
doesn't affect the weak type of the realization, and which guarantees that, when
the monoid $U$ is generated, it becomes a group.

Let $I$ be the category with two objects and one morphism $0\rightarrow 1$, and
let $\overline{I}$ be the category with two objects and an isomorphism between
them. Consider these as Segal categories (taking their nerve as
bisimplicial sets constant in the second variable).  Note that $|I|$ and $|
\overline{I}|$ are both contractible, so the obvious inclusion
$I\hookrightarrow \overline{I}$ induces an equivalence of realizations.

The bisimplicial set $\overline{I}$ is just that which is represented by
$(1,0)\in \Delta \times \Delta$. Thus for a Segal precat $A$,
if $f\in A_{1,0}$ is an object of $A_{1/}$ (a ``morphism'' in $A$) then it
corresponds to a morphism $I\rightarrow A$.
Set
$$
A^f:= A \cup ^{I} \overline{I}.
$$
Now the morphism $f$ is strictly invertible in the precategory $\tau _{\leq
1}(A^f)$ and in particular, when we apply the operations described above, the
image of $f$ becomes invertible in the resulting category.  If $A_0=\ast$
(whence $A^f_0=\ast$ too) then
the image of $f$ becomes invertible in the resulting monoid.
Note finally that
$$
|A | \rightarrow | A^f| = | A| \cup ^{|I|} |\overline{I}|
$$
is a weak equivalence.
In fact we want to invert all of the $1$-morphisms.
Let
$$
A' := A \cup ^{\bigcup _fI}\left( \bigcup _f \overline{I} \right)
$$
where the union is taken over all $f\in A_{1,0}$. Again
$|A|\rightarrow | A'|$ is a weak equivalence. Now, when we apply the previous
procedure to $\tau _{\leq 1} (A')$ giving a category $U$ (a monoid if $A$ had
only one object), all morphisms coming from $A_{1,0}$ become invertible. Note
that the morphisms in $A'$, i.e. objects of $A'_{1,0}$, are either morphisms in
$A$ or their newly-added inverses. Thus all of the morphisms coming from
$A'_{1,0}$ become invertible in the category $U$. But it is clear from the
operations described above that $U$ is generated by the morphisms in
$A'_{1,0}$. Therefore $U$ is a groupoid. In the case of only one object,
$U$ becomes a group.

By Segal's theorem we then have $U= \pi _1(| A|)$. If we know for some reason
that $|A|$ is simply connected, then $U$ is the trivial group. More
precisely, search for a proof that $\pi _1= 1$, and when such a proof is found,
apply the corresponding series of operations to
$\tau _{\leq 1} (A')$ to obtain $U=\ast$. Applying the
operations to $A'$ upstairs, we obtain a new $A''$ with $|A''|\cong |A'|\cong
|A|$ and  $A''_{1/}$ connected.

Another way of looking at this is to say that every time one needs to take the
inverse of an element in the proof that the group is trivial, add on a copy of
$\overline{I}$ over the corresponding copy of $I$.

\numero{An algorithm for calculating the $\pi _i(X)$, $X$ simply connected}

We describe how to use the above description of $\Omega X$ inductively to
obtain the $\pi _i(X)$.  This seems to be a new algorithm, different from those
of E. Brown \cite{EBrown} and Kan--Curtis \cite{Kan1} \cite{Kan2}
\cite{Curtis}.

There is an unboundedness to the resulting algorithm, coming essentially from a
problem with $\pi _1$ at each stage. Even though we know in advance that the
$\pi _1$ is abelian, we would need to know ``why'' it is abelian in a
precise way
in order to specify a strategy for making $A_{1/}$ connected at the appropriate
place in the loop.  In the absence of a particular description of the proof we
are forced to say ``search over all proofs'' at this stage.

\subnumero{The case of finite homotopy groups}

We first present our algorithm for the case of finite homotopy groups.
Suppose we want to calculate $\pi _n(X)$. We assume known that the
$\pi _i(X)$ are finite for $i\leq n$.

\noindent
{\em Start:}
Fix $n$ and
start with a simplicial set $X$ containing a finite number of nondegenerate
simplices. Suppose we know that $\pi
_1(X,x)$ is a given finite group; record this group, and set $Y$ equal to
the corresponding covering space of $X$. Thus $Y$ is simply connected. Now
contract out a maximal tree to obtain $Z$ with $Z_0=\ast$.

\noindent
{\em Step 1.} Let $A_{p,k}:= Z_p$ be the
corresponding Segal precat. It has only one object.

\noindent
{\em Step 2.} Let $A'$ be the coproduct of $A$ with one copy of the nerve of
the category $\overline{I}$ (containing two isomorphic objects), for each
morphism $I\rightarrow A$ (i.e. each point in $A_{1,0}$).

\noindent
{\em Step 3.} Apply the procedure of \S 3 to obtain a morphism $A'\rightarrow
A''$ with $A''_{1/}$ connected, and inducing a weak equivalence on
realizations.  (This step can only be bounded if we have a specific proof that
$\pi _1(Y,y)=1$).

\noindent
{\em Step 4.} Apply the procedure of \S 2 to obtain a morphism $A''
\rightarrow  B$ (inducing a weak equivalence on realizations) such that $B$ is a
Segal groupoid.
By the discussion of \S 2 and Theorem \ref{bound1},  the
$n-1$-type of $B_{1/}$ is effectively calculable.  By Segal's theorem,
$$
|B_{1/}| \sim \Omega | B| \sim \Omega | Y|,
$$
which in turn is the
connected component of $\Omega | X|$. Thus
$$
\pi _n(|X|)= \pi _{n-1}(|B_{1/}|).
$$

\noindent
{\em Step 5.} Go back to the {\em Start} with the new $n$ equal to the old
$n-1$, and the new $X$ equal to the simplicial set $B_{1/}$ above. The new
fundamental group is known to be abelian (since it is $\pi _2$ of the previous
$X$). Thus we can calculate the new fundamental group as $H_1(X)$ and, under our
hypothesis, it will be finite.

Keep repeating the procedure until we get down to $n=1$ and have recorded the
answer.

\subnumero{How to get rid of free abelian groups in $\pi _2$}

In the case where the higher homotopy groups are infinite (i.e. they contain
factors of the form $\zz ^a$) we need to do something to get past these
infinite groups. If we go down to the case where $\pi _1$ is infinite, then
taking the universal covering no longer results in a finite complex. We prefer
to avoid this by tackling the problem at the level of $\pi _2$, with a
geometrical argument. Namely, if $H^2(X, \zz)$ is nonzero then we can take a
class there as giving a line bundle, and take the total space of the
corresponding $S^1$-bundle. This amounts to taking the fiber of a map
$X\rightarrow K(\zz , 2)$.  This can be done explicitly and effectively,
resulting again in a calculable finite complex. In the new complex we will have
reduced the rank of $H_2(X, \zz )= \pi _2(X)$ (we are assuming that $X$ is
simply connected).

The original method of E. Brown \cite{EBrown} for effectively calculating the
$\pi _i$ was basically to do this at all $i$. The technical problems in
\cite{EBrown} are caused by the fact that one doesn't have a finite complex
representing $K(\zz , n)$.  In the case $n=2$ we don't have these technical
problems because we can look at circle fibrations and the circle is a finite
complex.  For this section, then, we are in some sense reverting to an easy
case of \cite{EBrown} and not using the Seifert-Van Kampen technique.

Suppose $X$ is a simplicial set with finitely many nondegenerate simplices,
and suppose $X_0=X_1 = \ast$.  We can calculate $H^2(X, \zz )$ as
the kernel of the differential
$$
d: \zz ^{X'_2} \rightarrow \zz ^{X'_3}.
$$
Here $X'_i$ is the set of nondegenerate $i$-simplices.
(Note that a basis of this kernel can effectively be computed using Gaussian
elimination). Pick an element $\beta$ of this basis, which is a collection of
integers $b_t$ for each $2$-simplex (i.e. triangle) $t$.  For each triangle $t$
define an $S^1$-bundle $L_t$ over $t$ together with trivialization
$$
L_t |_{\partial t} \cong \partial t \times S^1.
$$
To do this, take $L_t = t \times S^1$ but change the trivialization along the
boundary by a bundle automorphism
$$
\partial t \times S^1 \rightarrow \partial t \times S^1
$$
obtained from a map $\partial t \rightarrow S^1$ with winding number $b_t$.
 Let $L^{(2)}$ be the $S^1$-bundle over the $2$-skeleton of $X$ obtained by
glueing together the $L_t$ along the trivializations over their boundaries. We
can do this effectively and obtain $L^{(2)}$ as a simplicial set with a finite
number of nondegenerate simplices.

The fact that $d(\beta )=0$ means that for a $3$-simplex $e$, the restriction
of $L^{(2)}$ to $\partial e$ (which is topologically an $S^2$) is a trivial
$S^1$-bundle. Thus $L^{(2)}$ extends to an $S^1$-bundle $L^{(3)}$ on
the $3$-skeleton of $X$. Furthermore, it can be extended across any simplices
of dimension $\geq 4$ because all $S^1$-bundles on $S^k$ for $k\geq 3$, are
trivial ($H^2(S^k, \zz )= 0$). We obtain an $S^1$-bundle $L$ on $X$.
By subdividing things appropriately (including possibly subdividing $X$)
we can assume that $L$ is a simplicial set with a finite number of
nondegenerate simplices.  It depends on the choice of basis element $\beta$, so
call it $L(\beta )$.

Let
$$
T=L(\beta _1)\times _X \ldots \times _X L(\beta _r)
$$
where $\beta _1,\ldots , \beta _r$ are our basis elements found above.
It is a torus bundle with fiber $(S^1)^r$. The long exact homotopy sequence for
the map $T\rightarrow X$ gives
$$
\pi _i (T)= \pi _i (X), \;\;\; i \geq 3;
$$
and
$$
\pi _2(T) = \ker (\pi _2(X) \rightarrow \zz ^r).
$$
Note that $\zz ^r$ is the dual of $H^2(X , \zz )$ so the kernel $\pi _2(T)$ is
finite. Finally, $\pi _1(T)=0$ since the map  $\pi _2(X) \rightarrow \zz ^r$
is surjective.

Note that we have a proof that $\pi _1(T)=0$.

\subnumero{The general algorithm}

Here is the general situation. Fix $n$. Suppose $X$ is a simplicial set with
finitely many nondegenerate simplices, with $X_0=\ast$ and with a proof that $\pi
_1(X)=\{1\}$.  We will calculate $\pi _i(X)$ for $i\leq n$.

\noindent
{\em Step 1.}\,
Calculate (by Gaussian elimination) and record $\pi _2(X) = H_2(X, \zz )$.

\noindent
{\em Step 2.}\,
Apply the operation described in the previous subsection above, to obtain a new
$T$ with $\pi _i(T)=\pi _i(X)$ for $i\geq 3$, with $T_0 = \ast$, with $\pi
_1(T)=1$, and with $\pi _2(T)$ is finite.

Let $A$ be the  Segal precat corresponding to $T$.

\noindent
{\em Step 3.}\,
Use the discussion of \S 3 to obtain a morphism $A\rightarrow A'$ inducing a
weak equivalence of realizations, such that $A'_{1/}$ is connected. For this
step we need a proof that $\pi _1(T)=1$.  In the absence of a specific (finite)
proof, search over all proofs.

\noindent
{\em Step 4.}\,
Use Corollary \ref{thinkIcan} to replace $A'$ by a Segal  precat $B$ with
$|A'|\rightarrow |B|$ a weak equivalence, such that
the $n-1$-type of $B_{1/}$ is equivalent to $\Omega |B|$ which in turn is
equivalent to $\Omega |X|$. Let $Y= B_{1/}$ as a simplicial set.

Note that $Y$ is connected and $\pi _1(Y)$ is finite, being equal to
$\pi _2(T)$. We have $\pi _i(X) = \pi _{i-1}(Y)$ for
$3\leq i \leq n$.

\noindent
{\em Step 5.}\,
Choose a universal cover of $Y$, and mod out by a maximal
tree in the $1$-skeleton to obtain a simplicial set $Z$, with finitely many
nondegenerate simplices, with $Z_0 = \ast$, and with a proof that $\pi _1(Z)=
1$.  We have $\pi _i(X) = \pi _{i-1}(Z)$ for
$3\leq i \leq n$.

Go back to the beginning of the algorithm and plug in $(n-1)$
and $Z$. Keep doing this until, at the step where we calculate $\pi _2$ of the
new object, we end up having calculated $\pi _i(X)$ as desired.

\subnumero{Proofs of Godement}

We pose the following question: how could one obtain, in the process of
applying the above algorithm, an explicit proof that at each stage the
fundamental group (of the universal cover $Z$ in step 5) is trivial? This could
then be plugged into the machinery of \S 3 to obtain an explicit strategy,
thus we would avoid having to try all possible strategies. To do this we would
need an explicit proof that $\pi _1(Y)$ is finite in step 4, and this in turn
would be based on a proof that $\pi _1(Y)=\pi _2(T)$ as well as a proof of the
Godement property that $\pi _2(T)$ is abelian.

\numero{Example: $\pi _3(S^2)$}

The story behind this paper is that Ronnie Brown came by for Jean Pradines'
retirement party, and we were discussing Seifert-Van Kampen. He pointed out that
the result of \cite{nCAT} didn't seem to lead to any actual calculations.
After that, I tried to use that technique (in its simplified
Segal-categoric version) to calculate $\pi _3(S^2)$. It was apparent from this
calculation that the process was effective in general.

We describe here what happens for calculating $\pi _3(S^2)$.
We take as simplicial model a simplicial set
with the basepoint as unique $0$-cell $\ast$ and with one nondegenerate simplex
$e$ in degree $2$.  Note that this leads to many degenerate simplices in degrees
$\geq 2$ (however there is only one degenerate simplex which we denote $\ast$ in
degree $1$).

We follow out what happens in a language of cell-addition. Thus we don't feel
required to take the whole cone $C$ at each step of an operation $Arr (A,m)$;
we take any addition of cells to $A_{m/}$ lifting cells in $A_{1/}\times \ldots
\times A_{1/}$.

We keep the notation $A$ for the result of each operation (since our discussion
is linear, this shouldn't cause too  much confusion).

Finally, we will pre-suppose the technique of the proof of Theorem
\ref{bound1} in \S 6 below (alternatively, the reader may take the present
example as an introduction to the following section).

The first step is to $(2,0)$-arrange $A$. We do this by adding a $1$-cell
joining the two $0$-cells in $A_{2/}$, in an operation of type $Arr (A,
2)$. Note that both $0$-cells map to the same point $A_{1/}\times A_{1/} =
\ast$.
The first result of this is to add on $1$-cells in the $A_{m/}$ connecting all
of the various degeneracies of $e$, to the basepoint. Thus the $A_{m/}$ become
connected.
Additionally we get a new $1$-cell added onto to $A_{1/}$
corresponding to the third face $(02)$.  Furthermore, we obtain all images of
this cell by degeneracies $m\rightarrow 1$. Thus we get $m$ circles
attached to  the pieces which became connected in the first part of this
operation. Now each $A_{m/}$ is a wedge of $m$ circles.

In particular note that $A$ is now $(m,1)$-arranged for all $m$.

The next step is to $(2,2)$-arrange $A$. To do this, note that the Segal map
is
$$
S^1 \vee S^1 = A_{2/} \rightarrow A_{1/} \times A_{1/} = S^1 \times S^1.
$$
To arrange this map we have to add a $2$-cell to $S^1 \vee S^1$
with attaching map the commutator relation. Again, this has the result of
adding on $2$-cells to all of the $A_{m/}$ over the pairwise commutators of the
loops. Furthermore, we obtain an extra $2$-cell added onto $A_{1/}$ via the
edge $(02)$. The attaching map here is the commutator of the generator with
itself, so it is homotopically trivial and we have added on a $2$-sphere.
(Note in passing that this $2$-sphere is what gives rise to the class of the
Hopf map).
 Again, we obtain the  images of this $S^2$ by all of the degeneracy
maps $m\rightarrow 1$. Now
$$
A_{1/} = S^1 \vee S^2,
$$
$$
A_{2/} = (S^1 \times S^1)\vee S^2\vee S^2,
$$
and in general $A$ is $(m, 2)$-arranged for all $m$ (the reader is invited
to check for himself).  Looking forward to the next section, we see that adding
$3$-cells to $A_{m/}$ for $m\geq 3$  in the appropriate way as described in the
proof of \ref{bound1}, will end up resulting in the addition of $4$-cells (or
higher) to $A_{1/}$ so this no longer affects the $2$-type of $A_{1/}$. Thus
(for the purposes of getting $\pi _3(S^2)$) we may now ignore the $A_{m/}$ for
$m\geq 3$.

The remaining operation is to $(2,3)$-arrange $A$. For this, look at the
Segal map
$$
A_{2/} = (S^1 \times S^1)\vee S^2\vee S^2 \rightarrow
$$
$$
A_{1/}\times A_{1/} = (S^1 \vee S^2)\times (S^1 \vee S^2).
$$
Let $C$ be the mapping cone on this map. Then we end up attaching one copy
of $C$ to $A_{1/}$ along the third edge map $A_{2/} \rightarrow A_{1/}$.
This gives the answer for the $2$-type of $\Omega S^2$:
$$
\tau _{\leq 2}(\Omega S^2) =
\tau _{\leq 2}\left(
(S^1 \vee S^2)\cup ^{(S^1 \times S^1)\vee S^2\vee S^2}C \right) .
$$
To calculate $\pi _2(\Omega S^2)$ we revert to a homological formulation
(because it isn't easy to ``see'' the cone $C$).  In homology of degree $\leq
2$,  the above Segal map
$$
(S^1 \times S^1)\vee S^2\vee S^2 \rightarrow
(S^1 \vee S^2)\times (S^1 \vee S^2)
$$
is an isomorphism. Thus the map $A_{2/} \rightarrow C$ is an isomorphism on
homology in degrees $\leq 2$, and adding in a copy of $C$ along $A_{2/}$
doesn't change the homology. Thus
$$
H_2(\Omega S^2) = H_2 (S^1 \vee S^2) = \zz .
$$
Noting that (as we know from general principles) $\pi _1(\Omega S^2) = \zz
$ acts
trivially on $\pi _2(\Omega S^2)$ and $\pi _1$ itself has no homology in degree
$2$, we get that $\pi _2(\Omega S^2)= H_2(\Omega S^2) = \zz$.

{\em Exercise:} Calculate $\pi _4(S^3)$ using the above method.

{\em Remark:} our above recourse to homology calculations suggests that it
might be interesting to do pushouts and the operation $Cat$ in the context of
simplicial chain complexes.

\subnumero{Seeing Kan's simplicial free groups}
Using the above procedure, we can actually see how Kan's simplicial free groups
arise in the calculation for an arbitrary simplicial set $X$. They
arise just from a first stage where we add on $1$-cells. Namely, if in doing the
procedure $Arr (A,m)$ we replace $C$ by a choice of $1$-cell joining any two
components of $A_{m/}$ which go to the same component under the Segal map, then
applying this operation for various $m$, we obtain a simplicial space whose
components are connected and homotopic to wedges of circles.  (We have to start
with an $X$ having $X_1=\ast$). The resulting simplicial space has the same
realization as $X$. If $X$ has only finitely many nondegenerate simplices then
one can stop after a finite number of applications of this operation.  Taking
the fundamental groups of the component spaces (based at the degeneracy of the
unique basepoint) gives a simplicial free group. Taking the classifying
simplicial sets of these groups in each component we obtain a bisimplicial set
whose realization is equivalent to $X$. This bisimplicial set actually satisfies
$A_{p,0}= A_{0,k}=\ast$, in other words it satisfies the globular condition in
both directions!  We can therefore view it as  a Segal precat in two ways. The
second way, interchanging the two variables, yields a Segal precat where the
Segal maps are {\em isomorphisms} (because at each stage it was the classifying
simplicial set for a group).  Thus, viewed in this way, it is a Segal groupoid
and Segal's theorem implies that the simplicial set $p\mapsto A_{p,1}$, which
is the underlying set of a simplicial free group, has the homotopy type of
$\Omega X$.

\numero{Proof of Theorem \ref{bound1}}

Note that Theorem \ref{bound1} is just used to prove that our
process works. Thus we may use an infinite procedure in its proof.

Here is the idea. Let $\Gamma (m)\subset h(m)$ denote some ``horn''
, i.e. union of all faces but one.
\footnote{
Preferably
a horn which is compatible with the ordering of the product cf \cite{Street},
i.e. taking out some face except the first or last. This condition is so that
our arranging operation will be compatible with the notion of Segal
category---where some arrows may not be invertible---and not just Segal
groupoid.}
Use the notations of \S 1. For any inclusion of simplicial sets
$B'\subset B$ we can set $$ U:=(\Gamma (m) \otimes B) \cup ^{\Gamma (m) \otimes
B'} (h(m) \otimes B')
$$ and
$$
V:= h(m)\otimes B,
$$
we obtain an inclusion of bisimplicial sets $U\subset V$ such that
$|U|\rightarrow |V|$ is a weak equivalence. If $A$ is a Segal precat and
$U\rightarrow A$ a morphism then set $A':= A \cup ^UV$. The morphism
$|A|\rightarrow |A' |$ is a weak equivalence. Note that $A'$ is again a Segal
precat because the morphism $\Gamma (m)\rightarrow h(m)$ includes all of the
vertices of $h(m)$.  Finally, note that for $p\leq m-2$
the morphism
$$
A_{p/}\rightarrow A'_{p/}
$$
is an isomorphism (because the same is true of $\Gamma (m)\rightarrow h(m)$).
This last property allows us to conserve the homotopy type of the smaller
$A_{p/}$.

We would like to use operations of the above form, to $(m,k)$-arrange $A$.
In order to do this we analyze what a morphism from $U$ to $A$ means.

For any simplicial set $X$, we can form the simplicial set ${[}X, A]$ with the
property that a map $B\rightarrow {[}X,A]$ is the same thing as a map $X\otimes
B\rightarrow A$.  If $X\hookrightarrow Y$ is an inclusion obtained by adding
on an $m$-simplex $h(m)$ over a map $Z\rightarrow X$ where $Z\subset h(m)$ is
some subset of the boundary, then
$$
{[}Y,A] = {[}X , A] \times _{{[}Z,A]} A_{m/}.
$$
In this way, we can reduce ${[}X,A]$ to a gigantic iterated fiber product of the
various components $A_{p/}$.

{\em Claim:} there exists a cofibration $A\rightarrow A'$ such that $A_{p/}
\rightarrow A'_{p/}$ is a weak equivalence for all $p$, and such that $A'$ is
{\em fibrant} in the sense that for any cofibration of simplicial sets $X\subset
Y$ including all of the vertices of $Y$, the morphism
$$
{[}Y,A'] \rightarrow {[}X,A']
$$
is a Kan fibration of simplicial sets. In order to construct $A'$, we just
``throw in'' everything that is necessary. More precisely, suppose
$B'\subset B$ is a trivial fibration of simplicial sets. A diagram
$$
\begin{array}{ccc}
B' & \rightarrow & B \\
\downarrow && \downarrow \\
{[}Y,A'] &\rightarrow & {[}X,A']
\end{array}
$$
corresponds to a diagram
$$
\begin{array}{ccc}
U & \rightarrow & X \otimes B \\
\downarrow && \downarrow \\
A' &\rightarrow & A'
\end{array}
$$
with $U= (Y\otimes B')\cup ^{X\otimes B'} (X\otimes B)$.
The morphism of bisimplicial sets $U\rightarrow X \otimes B$ is a weak
equivalence on each vertical column $(\, )_{p/}$.  Therefore we can throw in
to $A'$ the pushout along this morphism, without changing the weak equivalence
type of the $A_{p/}$.  Note that the new thing is again a Segal precat
because of the assumption that $X\subset Y$ contains all of the vertices. Keep
doing this addition over all possible diagrams, an infinite number of times,
until we get the required Kan fibration condition to prove the claim. (I have
expressed this argument in an informal way. The best way to write it up
formally speaking is certainly in a closed model category formalism---which is
left for another time.)

We now describe a procedure for proving Theorem \ref{bound1}. At each step of
the procedure, we apply
(without necessarily saying so everywhere) the construction of the previous
claim to recover the  fibrant condition. Thus  we may always assume that
our Segal precat satisfies the fibrant condition of the previous paragraph.

We have already described above an arranging operation $A\mapsto Arr (A,m)$.
We now describe a second arranging operation, under the hypothesis that
$A$ is a fibrant (in the sense of the previous claim) Segal precat.
Fix $m$ and fix a horn $\Gamma (m)\subset h(m)$ (complement of all but one of
the faces, and the face that is left out should be neither the first nor the
last face).  Let $C$ be the cone on the map
$$
A_{m/} = {[}h(m), A] \rightarrow {[}\Gamma (m), A].
$$
Thus we have a diagram
$$
A_{m/} \rightarrow C \rightarrow {[}\Gamma (m), A].
$$
Note that $\Gamma (m)$ is a gigantic iterated fiber product of various $A_{p/}$
for $p<m$. This diagram corresponds to a map $U\rightarrow A$
where
$$
U:= (\Gamma (m) \otimes C )\cup ^{\Gamma (m)\otimes A_{m/}} h(m)\otimes A_{m/}.
$$
Letting $V:= h(m)\otimes C$ we set
$$
Arr2(A, m):= A\cup ^U V.
$$
Notice first of all that by the previous discussion, the map
$$
|A| \rightarrow | Arr2(A, m)|
$$
is a weak equivalence of spaces.

We have to try to figure out what effect
$Arr2(A,m)$ has.
We do this under the following hypothesis on the utilisation of this operation:
that $A$ is $(m, k-1)$-arranged, and $(p,k)$-arranged for all $p<m$.

The first step is to notice that the fiber product in the expression of
${[}\Gamma (m), A]$ is a homotopy fiber product, because of the ``fibrant''
condition we have imposed on $A$. Furthermore the elements in this fiber
product all satisfy the Segal condition up to $k$ (bijectivity on $\pi _i$
for $i<k$ and surjectivity for $\pi _k$). Thus the morphism
$$
{[}\Gamma (m) , A]
\rightarrow A_{1/} \times _{A_0} \ldots \times _{A_0} A_{1/}
$$
(the Segal fiber product for $m$, on the right)
is an isomorphism on $\pi _i$ for $i<k$ and a surjection for $i=k$.
Thus when we add to $A_{m/}$ the cone $C$, we obtain the condition of being
$(m,k)$-arranged.

By hypothesis, $A$ is $(m,k-1)$-arranged, in particular the map
$A_{m/} \rightarrow C$ is an isomorphism on $\pi _i$ for $i<k-1$ and
surjective for $\pi _{k-1}$.  Thus $C$ may be viewed as obtained from $A_{m/}$
by adding on cells of dimension $\geq k$.  Therefore, for all $p$ the morphisms
$U_{p/}\rightarrow V_{p/}$ are homotopically obtained by addition of cells of
dimension $\geq k$.

From the previous paragraph, some extra cells of dimension $\geq k$ are added to
various $A_{p/}$ in the process. This doesn't spoil the condition of being
$(p,k)$-arranged wherever it exists.  However (and here is the major advantage
of this second operation) the $A_{p/}$ are left unchanged for $p\leq m-2$.  This
is because all $p$-faces of the $m$-simplex are then contained in the horn
$\Gamma (m)$.

We review the above results. First, the hypotheses on $A$ were:
\newline
(a)\,\, that $A$ is fibrant in the sense of the claim above;
\newline
(b)\,\, that $A$ is $(p,k)$-arranged for $p<m$, and $(m,k-1)$-arranged.
We then obtain a construction $Arr2(A,m)$ with the following properties:
\newline
(1)\,\, the map
$A_{p/}\rightarrow Arr2(A,m)_{p/}$ is an isomorphism for $p\leq m-2$;
\newline
(2)\,\, for any $p$ the map $A_{p/}\rightarrow Arr2(A,m)_{p/}$ induces an
isomorphism on $\pi _i$ for $i<k-1$;
\newline
(3)\,\, if $A$ is $(p,k)$-arranged for any $p$ then $Arr2(A,m)$ is also
$(p,k)$-arranged;
\newline
(4)\,\, and $Arr2(A,m)$ is $(m,k)$-arranged.

{\em Remark:} at $m=2$ the operations $Arr(A,2)$ and $Arr2(A,2)$ coincide.

With an infinite series of applications of the construction $Arr2(A,m)$ and the
fibrant replacement operation we can prove Theorem \ref{bound1}. The reader may
do this as an exercise or else read the explanation below.

Take an array of dots, one for each $(p,k)$. Color the dots green if $A$ is
$(p,k)$-arranged, red otherwise (note that one red dot in a column implies
red dots everywhere above).  We do a sequence of operations of the form
fibrant replacement (which doesn't change anything) and then $Arr2(A,m)$.
When we do this, change the colors of the dots appropriately.

Also mark an
$\times$ at any dot $(p,k)$ such that the $\pi _i(A_{p/})$ change for
any $i \leq k-1$.  (Keep any $\times$ which are marked, from one step to
another).  If a dot $(p,k)$ is never marked with a $\times$ it means that
the $\pi _i(A_{p/})$ remain unchanged for $i<k$.

We don't color the dots $(1,k)$ but we still might mark an
$\times$.

Suppose the dot $(m,k)$ is red,
the dots $(p,k)$ are green for $p<m$ and the dot $(m,k-1)$ is green.
Then apply the fibrant replacement and the operation $Arr2(A	,m)$.
This has the following effects.  Any green dot $(p,j)$ for $p\leq m-2$
(and arbitrary $j$) remains green. The dot $(m-1,k)$ remains green.
However, the dot $(m-1, k+1)$ becomes red.
The dot $(m,k)$ becomes green.
The dots $(m-1, k)$ and $(m,k)$, as well as all $(p,k)$ for $p>m$, are marked
with an $\times$. The dots above these are also marked with an $\times$ but
no other dots are (newly)  marked with an $\times$.

In the situation of Theorem \ref{bound1}, we start with green dots at $(p,k)$
for $p+k \leq n$. We may as well assume that the rest of the dots are colored
red. Start with $(m,k)=(n+1, 0)$ and apply the procedure of the previous
paragraph. The dot $(n+1,0)$ becomes green, the dot $(n,0)$ stays green,
and the dots $(n, 0), (n+1, 0), \ldots $ are marked with an $\times$. Continue
now at $(n,1)$ and so on. At the end we have made all of the dots $(p,k)$ with
$p+k=n+1$ green, and we will have marked with an $\times$ all of the dots
$(p,k)$ with $p+k=n$ (including the dot $(1,n-1)$; and also all of the dots
above this line).

We can now iterate the procedure. We successively get green dots on each of the
lines $p+k= n+j$ for $j=1,2,3,\ldots$. Furthermore, no new dots will be marked
with a $\times$.  After taking the union over all of these iterations, we obtain
an $A'$ which is $(p,k)$-arranged for all $(p,k)$.
Thus $A'$ is a Segal category.

Note that the morphism $|A|\rightarrow |A'|$ is  a weak equivalence of spaces.

By looking at which dots are
marked with an $\times$, we find that the morphisms
$$
A_{p/} \rightarrow A'_{p/}
$$
induce isomorphisms on $\pi _i$ whenever $i < n-p$.   This completes the proof
of Theorem \ref{bound1}.
\eop

\numero{Complementary remarks}

What we are doing above is essentially applying generalized Seifert-Van Kampen.
To explain this, we remark first of all that the category $SPC$ of Segal
precats has a structure of closed model category which we shall now explain
(refering to \cite{nCAT} for the proofs).

The cofibrations are the injections of
bisimplicial sets.

In order to define the notion of weak
equivalence, we need an operation $SeCat$ which transforms a Segal precat $A$
into a Segal category $SeCat(A)$. This can be defined, for example, as the
direct
limit of  an infinite number of applications of the operation $Arr(A,m)$.
(One has to be careful about using the operation $Arr2(A,m)$ defined in \S 6
because we need to pay attention to the direction of arrows---this was why we
didn't allow the horns obtained by taking out the first or last faces).
For the present abstract purposes it will be more convenient to define $SeCat$
in the following way (this mimics the discussion in \cite{nCAT} and we don't
repeat the proofs). We say that a Segal precat $A$ is {\em marked} if for every
inclusion of finite complexes $B'\hookrightarrow B$ and every map
$$
u: \Sigma (m) \otimes B \rightarrow A
$$
together with compatible extension
$$
v: h (m) \otimes B' \rightarrow A,
$$
we are given an extension of these two to a map
$$
f(B'\subset B, u, v): h(m)\otimes B \rightarrow A.
$$
It is clear that if $A$ is marked then it is a Segal category. It is also
clear that for any Segal precat $A$ there is a universal morphism
$$
\eta : A \rightarrow SeCat(A)
$$
to a marked Segal precat (universal in the sense that
if $B$ is any marked Segal precat and $A\rightarrow
B$ a morphism, this extends uniquely to a morphism $SeCat(A)\rightarrow B$). The
resulting functor $SeCat$ with natural transformation is a {\em monad} on the
category of Segal precats:

\begin{lemma}
\mylabel{monad}
The universal property applied to the identity
$SeCat(A)\rightarrow SeCat(A)$ gives a natural map
$$
SeCat(SeCat(A))\rightarrow SeCat(A)
$$
whose composition with either of the standard inclusions
$$
SeCat(A)\tworightarrows SeCat(SeCat(A))
$$
is the identity. Thus $SeCat$ is a {\em monad} on the category of Segal precats
(cf \cite{May}).

If $A$ is already a Segal category then the
morphism
$$
A\rightarrow SeCat(A)
$$
is an equivalence.
\end{lemma}
{\em Proof:}
See \cite{nCAT}, Lemma 3.4.
\eop

These properties characterize
the construction $SeCat$ up to equivalence (cf \cite{nCAT} Proposition 4.2,
but that must be a well-known technique).
The construction $SeCat(A)$ can be expressed using variants of the operations
$Arr (A,m)$ (cf \cite{nCAT}, \S 4).

We return to the discussion of the closed model structure.
Say that a morphism $f: A\rightarrow B$ is a {\em weak equivalence}
if $SeCat(A)\rightarrow SeCat(B)$ is an equivalence of Segal categories.
Finally
a morphism is said to be a {\em fibration} if it satisfies the lifting property
with respect to all cofibrations which are weak equivalences.

\noindent
{\em Caution:} the resulting notion
of ``fibrant'' is {\em not} the same as the notion refered to in \S 6 above
(the notion used there was just the convenient thing for the moment). The
``official'' notion of fibrant for Segal precats is the one from the present
paragraph.

\begin{theorem}
\mylabel{cmc}
The category of Segal precats with the above classes of cofibrations, weak
equivalences and fibrations, is a closed model category (and is ``internal'' in
the sense of \cite{nCAT}).
\end{theorem}

The proof of this is
essentially the same as that given in \cite{nCAT} and we don't repeat it here.

One can similarly define notions of Segal $2$-category and so on, and  at each
stage we obtain a  closed model category.

As happens in \cite{nCAT}, we obtain a Segal $2$-category $SECAT$ whose objects
are the fibrant Segal categories.

A similar argument to that of \cite{limits} shows that the fibrant replacement
$SECAT'$ of the Segal $2$-category of Segal categories, admits limits. In
particular it admits pushouts.

Here is how to calculate the pushout (similar to the calculation at the end
of \cite{limits}): if
$$
A \leftarrow B \rightarrow C
$$
is a diagram in $SECAT$ then it can be replaced by an equivalent diagram
(we use the same notation) where the morphisms are cofibrations of
fibrant Segal precats.  Note also that this can be done even if we start with
a diagram of Segal precats (technically speaking the objects of $SECAT$ are
fibrant Segal precats, but for calculational reasons we would like to look at
any Segal precat, identifying it with its fibrant replacement).

Assuming that the morphisms in the above diagram are cofibrations, let
$P:= A \cup ^B C$ be the pushout of Segal precats (i.e. the pushout of
bisimplicial sets).  To obtain the pushout in $SECAT$, replace $P$ by a weak
equivalent fibrant $P'$. By the definition of weak equivalence,
$$
SeCat (P) \rightarrow SeCat(P')
$$
is an equivalence of Segal categories. Note that $P'$ is already a Segal
category and $P' \rightarrow SeCat(P')$ is an equivalence of Segal categories
(the result analogous to this, for $n$-categories, is proved in \cite{nCAT}).
The conclusion of all of this is that to get a Segal category equivalent to
the pushout in $SECAT$, we just take
$$
SeCat(P)= SeCat(A\cup ^B C).
$$
Thus we may consider $SeCat(A\cup ^BC)$ as the {\em pushout of Segal
categories $ A$, $B$, and $C$}.

We say that a Segal category $A$ is {\em $1$-connected} if for all pairs of
objects $(x,y)$, $A_{1/}(x,y)$ is nonempty and connected.

The main result of this paper may be rewritten as follows.

\begin{theorem}
\mylabel{rewrite}
Suppose $A\leftarrow B \rightarrow C$ is a diagram of fibrations of Segal
groupoids such that $A$, $B$ and $C$ have only one object, and are
$1$-connected. Then the pushout $SeCat(A\cup ^BC)$ can be effectively calculated
within any finite region of $(m,k)$ (using the previous type of notation).
\end{theorem}

This situation is much better than the situation for the fundamental group,
where the specification of a group by generators and relations is known not to
be effective.

Now suppose $X$ is a topological space. Tamsamani defines a simplicial space
$\underline{\Omega} (X)$. If we apply to each component the singular simplicial
complex then we obtain a Segal precat which we denote by $\Pi _{\rm Seg}(X)$.
Thus denoting by $R^k = | h(k)|$ the standard topological $k$-simplex,
$$
\Pi _{\rm Seg}(X) _{p,k}= Hom ' (R^p\times R^k, X)
$$
where $Hom '$ denotes the subset of morphisms which are constant on $\{ v\}
\times R^k$ for all of the vertices $v$ of $X$.

This is already a Segal category, in fact a Segal groupoid (Tamsamani proves
that  $\underline{\Omega} (X)$ satisfies the Segal condition \cite{Tamsamani}).
We call it the {\em Poincar\'e fundamental Segal groupoid of $X$}.
Note that
$$
| \Pi _{\rm Seg}(X)_{1/}| \cong \Omega X
$$
is equivalent to the loop space of $X$.
Since the notion of Segal groupoid is expressly intended to be a delooping
machine, this is just the statement saying how a loop space has the delooping
structure.

If $X$ has a basepoint $x$ then we could define a based version $\Pi _{\rm
Seg}(X, x)$ by requiring that the vertices $\{ v\}
\times R^k$ get mapped to $x$. This gives a Segal groupoid with only one object.

If $X$ is simply connected then for any $n$
we can calculate $\Pi _{\rm Seg}(X,x)$
(up to equivalence of Segal groupoids)
within the region of dots $(m,k)$  for $m+k\leq n$. This is exactly the
procedure described above.

We have the following generalized Seifert-Van Kampen theorem for $\Pi _{\rm
Seg}(X)$.

\begin{theorem}
\mylabel{VK}
Suppose $U, V \subset X$ are open sets with $X=U\cup V$ and $W:=U\cap V$. Then
the diagram
$$
\begin{array}{ccc}
\Pi _{\rm Seg}(W) & \rightarrow & \Pi _{\rm Seg}(U) \\
\downarrow & & \downarrow \\
\Pi _{\rm Seg}(V) & \rightarrow & \Pi _{\rm Seg}(X)
\end{array}
$$
is a pushout of Segal groupoids. If $U$, $V$ and $W$ are connected then we may
choose a common basepoint and the diagram of based fundamental Segal groupoids
is a pushout. If, furthermore, $X$ is simply connected then the pushout can
effectively be computed.
\end{theorem}
{\em Proof:}
The realization of bisimplicial sets gives an essential inverse to $X\mapsto
\Pi _{\rm Seg}(X)$. This is an immediate consequence of Segal's Theorem
\ref{segal} (cf the argument in\cite{Tamsamani}).  Furthermore,
$\Pi _{\rm Seg}$ transforms injections of spaces to injections (cofibrations) of
Segal precats. Thus the upper and left vertical arrows of the above diagram are
cofibrations. Let
$$
\begin{array}{ccc}
\Pi _{\rm Seg}(W) & \rightarrow & \Pi _{\rm Seg}(U) \\
\downarrow & & \downarrow \\
\Pi _{\rm Seg}(V) & \rightarrow & P
\end{array}
$$
be the pushout of Segal precats (i.e. the pushout of bisimplicial sets).
Then
$$
\begin{array}{ccc}
|\Pi _{\rm Seg}(W)| & \rightarrow & |\Pi _{\rm Seg}(U)| \\
\downarrow & & \downarrow \\
| \Pi _{\rm Seg}(V)| & \rightarrow & |P|
\end{array}
$$
is a pushout of spaces (realization transforms pushouts of bisimplicial sets to
pushouts of spaces).
Thus $|P|$ is equivalent to $X$. We have a diagram of Segal
precats $$ \begin{array}{ccc}
P & \rightarrow & \Pi _{\rm Seg}(X) \\
\downarrow & & \downarrow \\
SeCat(P) & \rightarrow & SeCat(\Pi _{\rm Seg}(X)).
\end{array}
$$
The right vertical arrow is an equivalence of Segal categories (cf
\cite{nCAT}
Lemma 3.4).
In general the morphism $|A|\rightarrow |SeCat(A)|$ is a weak
equivalence of spaces. Therefore, in the induced diagram of realizations
$$
\begin{array}{ccc}
|P| & \rightarrow & |\Pi _{\rm Seg}(X)| \\
\downarrow & & \downarrow \\
|SeCat(P)|& \rightarrow & |SeCat(\Pi _{\rm Seg}(X))|
\end{array}
$$
the vertical arrows are weak equivalences. On the other hand, both
$|P|$ and $|\Pi _{\rm Seg}(X)|$ are equivalent to $X$ (by maps compatible with
the  upper arrow) so the upper arrow is a weak equivalence. This implies that
$$
|SeCat(P)| \rightarrow  |SeCat(\Pi _{\rm Seg}(X))|
$$
is a weak equivalence of spaces. However, since both
$SeCat(P)$ and $SeCat(\Pi _{\rm Seg}(X))$ are Segal groupoids, we have
$$
|SeCat(P)_{1/}(x,y)|\sim Path ^{x,y} |P|
$$
and
$$
|SeCat(\Pi _{\rm Seg}(X))_{1/}(x,y)|\sim Path ^{x,y} |\Pi _{\rm Seg}(X)|.
$$
This implies that the morphism
$$
SeCat(P) \rightarrow SeCat(\Pi _{\rm Seg}(X))
$$
is a fully faithful morphism of Segal groupoids.
The condition that $X= U\cup V$ means that the morphism is surjective on
objects. Therefore it is an equivalence, so
$SeCat(\Pi _{\rm Seg}(X))$ is a pushout of our original diagram.
Combined with the fact noted above
that
$$
\Pi _{\rm Seg}(X) \rightarrow SeCat(\Pi _{\rm Seg}(X))
$$
is an equivalence, we obtain that $\Pi _{\rm Seg}(X)$ is a pushout of our
original diagram. This completes the proof  (the proof in the pointed connected
case is the same).
\eop

We could say, use this theorem to calculate $\Pi _{\rm Seg}(X)$. We have
essentially done that, except that instead of taking pushout of various
contractible Segal groupoids one for each simplex, we just take pushout of  the
weakly contractible Segal precats $h(m)\otimes \ast$, which amounts to
saying, look at a simplicial set $X$ as a Segal precat constant in the second
variable. If we want to end up with a Segal groupoid we should start with a case
where the $1$-morphisms will be invertible.

More generally, starting with a simplicial set $X$ we can consider it as a
Segal precat constant in the second variable, and take $SeCat(X)$. This will
not in general be a Segal groupoid but only a Segal category: there is no
reason for the morphisms to be invertible.  If $X_0=\ast$ then it is a Segal
monoid.  We could call it the {\em Segal monoid generated by generators and
relations $X$}.

This process commutes with the classifying space construction. We have the
following result: if
$$
A \leftarrow B \rightarrow C
$$
is a diagram of Segal categories, then (changing things by an equivalence)
we can
assume that one of the maps is a cofibration. We have
$$
| A| \cup ^{| B|} | C| = | A\cup ^B C|\cong | Cat (A \cup ^B C)|.
$$
Thus the classifying space construction (realization) commutes with
pushout $Cat (A \cup ^B C)$. This way we get around the flatness problem
pointed out by Fiedorowicz (\cite{Fiedorowicz}, \S 4), for pushout of
topological
monoids.

\subnumero{Relationship with other delooping machines}

Segal's machine is only one of many ``machines'' describing the structure
necessary to deloop a space. The other main family of machines is known as the
May family, based on the notion of {\em operad}.  A fundamental result in the
theory was the passage between these two types of machines, see
\cite{MayThomason} \cite{Thomason}---and for a recent addition to the subject,
\cite{Dunn}. We briefly recall the technique used to pass between these
machines, as well as some basic formalism surrounding the May family of
machines.

An {\em operad} $\Cc$ is a collection of spaces $\Cc (j)$ together with a
``function replacement'' operation (and a few others) modelled on the structure
of the spaces $\Ee _X(j)$ of maps $X\times \ldots \times X \rightarrow X$.
See \cite{May} for the details.
A {\em $\Cc$-space} is an operad morphism $\Cc \rightarrow \Ee _X$ which is
thought of as an ``action'' of $\Cc$ on $X$.

(For our $1$-delooping purposes, we apparently don't need to consider the
symmetric group action, cf the remark in \cite{May} p. 27.)

There is a weaker version of the notion of $\Cc$-space which is fundamental in
the comparison theorems of \cite{MayThomason}, \cite{Thomason}. To an operad
$\Cc$ one associates the {\em category of operators} $\widehat{\Cc}$ which is a
topological category with a functor $\widehat{\Cc} \rightarrow \Delta ^o$
(cf \cite{MayThomason} \cite{Thomason} or \cite{Dunn} for the definitions).
Furthermore if $\Pi$ denotes the category used in the above references
(essentially generated by the principal face maps) then there is a lifting
$$
\Pi \rightarrow \widehat{\Cc} \rightarrow \Delta ^o.
$$

A special case is when $\Cc = \Mm$ is the operad with exactly one $j$-ary
operation for each $j$ (thought of as a $j$-fold composition). Then
$\widehat{\Mm} = \Delta ^o$. This yields back Segal's machine.

One now can define a {\em $\widehat{\Cc}$-precat}, as being a continuous functor
$A: \widehat{\Cc}\rightarrow Top$, such that the object $0$ goes to a
discrete set $A_0$, which we think of as the set of objects.
\footnote{
We ignore the question of {\em properness} in the sense of \cite{May},
\cite{MayThomason} \cite{Thomason} since this can be avoided by replacing $Top$
by the category of simplicial sets.
We keep the topological language here because the references for delooping are
written in that framework.}
As before denote the
space image of the object $p$ by $A_{p/}$.
The restriction of $A$ to $\Pi $ allows one to look at the {\em Segal maps}
(similarly to above). We say that $A$ is {\em
special} (in the usual terminology of
\cite{SegalTopology} \cite{MayThomason} \cite{Thomason} \cite{Dunn}) or a {\em
$\widehat{\Cc}$-category} (a more suggestive terminology) if the Segal maps are
equivalences.

The notion of $\Cc$-space used by May in \cite{May} is recovered as a
$\widehat{\Cc}$-precat (with only one object) such that the Segal maps are
{\em isomorphisms} (not just equivalences).

In case $\Cc = \Mm$ and $\widehat{\Cc} = \Delta ^o$ we recover the previous
notions of Segal precats and Segal categories.

Recall \cite{May} that an {\em $A_{\infty}$-operad} is an operad $\Cc$ such that
the functor $\widehat{\Cc} \rightarrow \Delta ^o$ induces an equivalence on
morphism spaces. For example, Stasheff's notion of
$A_{\infty}$-space is the notion of a $\Kk$-space for an appropriate
$A_{\infty}$-operad $\Kk$.

For any
$A_{\infty}$-operad $\Cc$ we obtain the analogues of the  notions which we
have discussed above for Segal categories.
(From here on we assume without further making this explicit that our operads
are $A_{\infty}$-operads).

For example, if $A$ is a
$\widehat{\Cc}$-category then we obtain
a category $\tau _{\leq 1} (A)$ whose nerve is
the simplicial set $p\mapsto \pi _0(A_{p/})$ (the structure of functor on
$\Delta^o$ is assured by the contractibility of the morphism
$\widehat{\Cc} \rightarrow \Delta ^o$---and the fact that this is the nerve of a
category is a consequence of the speciality condition).

As before we can define the simplicial set $\tau _{\leq 1} (A)$
even if $A$ is only a $\widehat{\Cc}$-precat.

Suppose $A$ is a $\widehat{\Cc}$-precat.
If $x_0,\ldots , x_p\in A_0$ then we denote by
$A_{p/}(x_0,\ldots , x_p)$ the inverse image of $(x_0,\ldots , x_p)\in
A_0^{p+1}$
by the morphism (which is well-defined since $A_0$ is a discrete set)
$A_{p/}\rightarrow A_0^{p+1}$. In particular, the space $A_{1/}(x,y)$ is
thought of as the space of morphisms from $x$ to $y$.

We say that a morphism $f:A\rightarrow B$ of $\widehat{\Cc}$-categories
(which for now just means a natural transformation of continuous functors
on $\widehat{\Cc}$) is {\em fully faithful} if for any pair of objects
$x,y\in A_0$ the morphism $A_{1/}(x,y)\rightarrow B_{1/}(f(x), f(y))$ is a
weak equivalence of spaces. We say that $f$ is {\em fully faithful} if the
induced morphism of categories
$$
\tau _{\leq 1} A \rightarrow  \tau _{\leq 1} B
$$
is surjective on isomorphism classes of objects. Finally, we say that a
morphism $f: A\rightarrow B$ is an {\em
equivalence of  $\widehat{\Cc}$-categories} (or just {\em equivalence})
if it is fully faithful and essentially surjective.

\begin{conjecture}
\mylabel{con1}
For any $A_{\infty}$-operad $\Cc$, there is an operation $Cat _{\Cc}$ going
from $\widehat{\Cc}$-precats to $\widehat{\Cc}$-categories such that
$Cat _{\Cc}$ is
a monad on the category of
$\widehat{\Cc}$-precats; and if $A$ is a $\widehat{\Cc}$-category then the
morphism $A \rightarrow  Cat _{\Cc}(A)$ is an equivalence of
$\widehat{\Cc}$-categories.  If $A_{1/}$ is connected then calculation of $Cat
_{\Cc}(A)$ is effective.
\end{conjecture}

The definition of the operation $Cat_{\Cc}$ should be similar to our
constructive discussion of the operation $SeCat$ above (and thus similar to the
operation $Cat$ of \cite{nCAT}).

{\bf Pushout:}
It is easy to define the pushout of $\widehat{\Cc}$-precats.  If
$$
A\leftarrow B \rightarrow C
$$
are two morphisms of $\widehat{\Cc}$-precats with the first one a {\em
cofibration} (i.e.
the morphisms $B_{p/}\rightarrow A_{p/}$ are cofibrations of spaces)
then the rule
$$
p\mapsto A_{p/} \cup ^{B_{p/}} C_{p/}
$$
defines a new $\widehat{\Cc}$-precat (the action of $\widehat{\Cc}$ being
defined in the obvious way) which we denote by $A\cup ^BC$.

If we assume Conjecture \ref{con1} then we can define the {\em pushout of
$\widehat{\Cc}$-categories}. If
$$
A\leftarrow B \rightarrow C
$$
are two morphisms of
$\widehat{\Cc}$-categories with the first morphism a cofibration, define their
{\em pushout} to be
$$
Cat_{\Cc}(A\cup ^BC).
$$

In the case being considered in \cite{May} and elsewhere, the
categories have only one object. In that case, we can revert to the
standard terminology and speak of {\em $\widehat{\Cc}$-spaces}.
Furthermore, they
are often supposed to be {\em grouplike} which just means that $\tau _{\leq
1}(A)$ is a groupoid (i.e. a group in the $1$-object case).  Pushout will
preserve these properties.

One might call a grouplike $\widehat{\Cc}$-category a {\em
$\widehat{\Cc}$-groupoid}.

\begin{conjecture}
\mylabel{mayVK}
The functor which associates to every space $X$ its loop space considered as a
grouplike $\widehat{\Cc}$-space, takes pushouts of connected spaces to pushouts
of $\widehat{\Cc}$-spaces, and more generally pushouts of (not necessarily
connected) spaces to pushouts of $\widehat{\Cc}$-groupoids.
\end{conjecture}

Finally, we make the following

\begin{conjecture}
\mylabel{mayCMC}
For a wide range of $A_{\infty}$-operads $\Cc$, there is a closed model
structure on the category of $\widehat{\Cc}$-precats, where the
cofibrations are as defined above, the weak equivalences are the morphisms
$A\rightarrow B$ such that the induced $Cat_{\Cc}(A)\rightarrow Cat_{\Cc}(B)$
is an equivalence of $\widehat{\Cc}$-categories, and where the fibrations are
the maps satisfying the lifting property with respect to trivial cofibrations.
\end{conjecture}

For $\widehat{\Mm}$
Conjectures \ref{con1}, \ref{mayVK} and \ref{mayCMC}
are just Lemma \ref{monad}, and Theorems \ref{rewrite}, \ref{cmc}  and \ref{VK}.

{\em Remark:}
The notion of pushout does actually appear in a certain example in \cite{May}
and \cite{CohenLadaMay}. Namely, there one considers the {\em free $\Cc$-space
generated by a given space $X$}.  In our present terms this means taking the
bisimplicial set $h(1)\otimes X$ and identifying the two endpoints to the
basepoint. The realization of this bisimplicial set is just the suspension
$\Sigma X$. On the other hand, applying the operation $Cat_{\Cc}$ should give
the free $\Cc$-space generated by $X$. This would prove that this free
$\Cc$-space is the loop space of $\Sigma X$---a fact that is underlying all of
\cite{May} and \cite{CohenLadaMay}, where this is used to calculate with $\Omega
\Sigma X$. This principle goes back to the work of James \cite{James}.

In contrast to the above notion of $\widehat{\Cc}$-category, which is
inherently ``weak'' and hence shouldn't depend on the choice of $\Cc$ up to
equivalence, the original notion of May \cite{May} (modified by relaxing the
requirement that there be only one object) is that of {\em $\Cc$-category} in
which we require the Segal maps to be isomorphisms.

\begin{conjecture}
\mylabel{con2}
For an appropriately chosen $A_{\infty}$-operad $\Cc$ (one conjectures that the
``little $1$-cube'' operad $\Cc _1$ of \cite{BoardmanVogt} \cite{May} should
work) in which the spaces fit together ``freely enough'', there is a
notion of pushout of $\Cc$-categories such that the loop space functor from
spaces to $\Cc$-groupoids takes pushouts to pushouts.
\end{conjecture}

It might be interesting to further weaken the notion of operad itself, to the
case where $\widehat{\Cc}$ is no longer a topological category but only, say, a
Segal category, or, why not, a $\widehat{\Cc}$-category! The notion of
$\widehat{\Cc}$-precat should also be replaced by a weak notion of functor
from $\widehat{\Cc}$ to $Top$.  This
last version poses some obvious circularity problems, and it is not clear
whether they can be resolved.
This seems to require a large amount of
effort and it is not clear whether it has any payoffs.

The above remarks (eventually including the preceding paragraph but one
would hope to avoid that) should point the way for an extension of the theory of
$n$-categories of \cite{Tamsamani} based on Segal's delooping machine, to
theories based on other types of $1$-delooping machines.  This might eventually
lead to comparisons with the other operad-based approaches to $n$-categories
\cite{BaezDolan} \cite{Batanin}.

\subnumero{Relation with the Poincar\'e fundamental $n$-groupoid $\Pi _n$}

We would like to be able to do something similar to what is done above, to
obtain effectively an $n$-groupoid $A$ equivalent to Tamsamani's fundamental
$n$-groupoid $\Pi _n(X)$, with the
part of $A$ within a certain region, finitely calculable. We would like to do
this, for example, under the hypothesis that the homotopy groups are finite.

Similarly we would like effectively to be able to calculate the pushout of
$n$-groupoids. Of course there is the usual type of problem with $\pi _1$,
which is even more difficult to avoid because, since there may be noninvertible
arrows, there doesn't seem to be an analogue of the universal cover which would
allow us to treat the problem in the finite case.

The
temptation is just to iterate our algorithm for constructing $\Omega X$,
applying
it again to each of the simplicial sets $A_{p/}$ and so on. The problem is that
in order to do this, we need a {\em functorial} construction. In the above,
every
time a nontrivial $\pi _1$ is encountered (necessarily abelian since coming from
the higher $\pi _i(X)$) one takes the universal cover to get back to a simply
connected case. But there is no functorial choice of universal cover, so some
more work must be done to make the calculation functorial.

In this connection it should be noted that when the $\pi _i (X)$ are finite for
$i\leq n$, Ellis' construction \cite{Ellis} gives a simplicial finite group with
the same $n$-type as $X$. Since simplicial groups are automatically Kan, we can
apply the simplicial version of Tamsamani's $\Pi _n$ to obtain an $n$-category
$A$ (it is even fibrant!) such that the components $A_M$ are finite sets.
This remark is the analogue for
$n$-categories of the remark made by Ellis ({\em op cit}) for $cat^n$-groups.
Thus we know by another method that it is possible to get a finite $\Pi _n(X)$.
The question we are asking here is whether the algorithm presented above
can also be made to do the same thing.

\end{document}